%% file: bourba.tex
\documentclass{bourbaki}

\date{Novembre 2000}
\bbkannee{53\` eme ann\'ee, 2000-2001}
\bbknumero{880}
\title{Factorisation faible des
applications birationnelles}
\subtitle{[d'apr\`es Abramovich, Karu, Matsuki,
W{\l}odarczyk et Morelli]}
\author{Laurent BONAVERO}
\address{Universit\'e Grenoble I\\
Institut Fourier\\
Lab. de Math\'ematiques Pures\\
URA 5582 du CNRS\\
F-38402 SAINT-MARTIN D'H\`ERES Cedex}
\email{bonavero@ujf-grenoble.fr}

\input xy
\xyoption{all}
\usepackage[all]{xy}
\usepackage{epsfig}

\newcommand{\ZZ}{{\mathbf Z}}

\newcommand{\QQ}{{\mathbf Q}}

\newcommand{\PP}{{\mathbf P}}
\newcommand{\KK}{{\mathbf K}}

\newcommand{\Ho}{\operatorname{Hom}}

\newcommand{\dra}{\dashrightarrow}
\newcommand{\Supp}{\operatorname{Supp}}
\newcommand{\Stab}{\operatorname{Stab}}

\newcommand{\pimul}{\operatorname{\pi-mult}}
\newcommand{\pimp}{\operatorname{\pi-mp}}
\newcommand{\para}{\operatorname{par}}
\newcommand{\st}{\operatorname{Star}}

\def\finpreuve
{\hskip 3pt \vrule height6pt width6pt depth 0pt}

\begin{document}

\maketitle

\section{Introduction}

Dans tout ce texte, $\KK$ d\'esigne un corps alg\'ebriquement
clos de caract\'eristique nulle (la non nullit\'e
de la caract\'eristique est essentielle dans les parties \ref{cobord}
et \ref{desing},
les constructions des parties \ref{morelli}
et \ref{struc toroid} sont cependant
valables en caract\'eristique arbitraire).

\subsection{Rappels sur les \'eclatements}

Si $n\geq 2$ et si $V$ est un $\KK$-espace vectoriel
de dimension $n$,  
on note $\PP_{\KK} (V)$ l'espace projectif
des droites vectorielles de $V$. Si $V =\KK ^n$,
on note simplement $\PP ^{n-1} = \PP _{\KK}(V)$.
Consid\'erons 
$$ B_0(\KK ^n) = \{ (v,d) \in  \KK ^n \times \PP ^{n-1} \, | \, v\in d\}.$$
Si $(x_1,\ldots,x_n)$ sont les coordonn\'ees sur $\KK^n$ et
$[y_1:\ldots:y_n]$ les coordonn\'ees homog\`enes associ\'ees
sur $\PP ^{n-1}$,
alors
$$ B_0(\KK ^n) = \{ (x,y)\in \KK ^n\times \PP ^{n-1}\, | \,\forall i,j \, ,\,
x_iy_j=x_jy_i \}.$$
Il en d\'ecoule que $B_0(\KK ^n)$ 
est une sous-vari\'et\'e alg\'ebrique ferm\'ee lisse de
$\KK ^n\times \PP ^{n-1}$.
La premi\`ere projection $\pi: B_0(\KK ^n)\to \KK^n$
est une application r\'eguli\`ere birationnelle
qui se restreint en un isomorphisme
$\pi:  B_0(\KK ^n) \setminus \pi ^{-1}(0)
\to \KK^n \setminus \{0\}$ avec 
$\pi ^{-1}(0)= \{0\}\times \PP ^{n-1}$.
L'application birationnelle $\pi: B_0(\KK ^n)\to \KK^n$
s'appelle {\em l'\'eclatement de $\KK^n$ en $0$} ;
$0$ est le {\em centre} de $\pi$ et $\pi ^{-1}(0)$
est le {\em diviseur exceptionnel} de $\pi$.

Plus g\'en\'eralement, si $Y$ est une sous-vari\'et\'e ferm\'ee
lisse d'une vari\'et\'e alg\'ebrique lisse $X$, il
y a une vari\'et\'e alg\'ebrique lisse $B_Y(X)$ 
et une application r\'eguli\`ere birationnelle
$\pi : B_Y(X) \to X$ qui se restreint en un 
isomorphisme
$\pi:  B_Y(X) \setminus \pi ^{-1}(Y)
\to X \setminus Y$ et 
$\pi ^{-1}(Y) \simeq \PP (N_{Y/X})$ o\`u $N_{Y/X}$
d\'esigne le fibr\'e normal de $Y$ dans $X$. 
L'application birationnelle $\pi: B_Y(X) \to X$
s'appelle {\em l'\'eclatement de $X$ le long de $Y$, ou de centre $Y$}
et $\pi ^{-1}(Y)$ est le {\em diviseur exceptionnel} de $\pi$. 
Moralement, on remplace chaque point $y$ de $Y$ par
l'espace projectif des directions normales \`a $Y$ dans $X$
passant par $y$.
Si la donn\'ee initiale est $B_Y(X)$, on dit encore
que $\pi: B_Y(X) \to X$ est une {\em contraction}
de centre $Y$ (\'eclatement et contraction d\'esignent donc
dans ce texte la m\^eme application birationnelle, ils
sont utilis\'es respectivement suivant que la donn\'ee initiale
est $X$ (resp. $B_Y(X)$) et la donn\'ee finale $B_Y(X)$
(resp. $X$)).    

\subsection{Enonc\'e du th\'eor\`eme principal}

C'est un probl\`eme classique depuis une trentaine d'ann\'ees 
de savoir si on peut d\'ecomposer une application birationnelle
entre deux vari\'et\'es alg\'ebriques compl\`etes 
lisses en une suite d'\'ecla\-te\-ments
et contractions de centres lisses. En dimension $1$,
la question est vide : toute application
birationnelle
entre deux courbes alg\'ebriques compl\`etes lisses est un isomorphisme. 
Dans le cas des surfaces, on sait depuis un si\`ecle
que toute application birationnelle entre surfaces compl\`etes lisses
est une suite d'\'eclatements et contractions de centre
des points (voir par exemple \cite{BPV84}).
Le probl\`eme \'etait ouvert d\`es la dimension trois. 

Nous donnons dans cet expos\'e
les grandes lignes de la d\'emonstration du th\'eor\`eme suivant,
d\^u \`a Abramovich, Karu, Matsuki,
W{\l}odarczyk \cite{AKMW99} et W{\l}odarczyk \cite{Wlod99}.

\begin{theo}\label{faible general}
Soit $\varphi:X_1 \dashrightarrow X_2$ 
une application birationnelle entre deux vari\'et\'es 
alg\'e\-bri\-ques compl\`etes et lisses
$X_1$ et $X_2$ sur $\KK$. 
Alors, $\varphi$ se factorise en 
une suite d'\'eclatements et de contractions de centres lisses.
Autrement dit, il y a une suite 
d'ap\-pli\-ca\-tions birationnelles entre vari\'et\'es alg\'ebriques
compl\`etes et lisses
$$ X_1 = V_0 \stackrel{\varphi_0}{\dashrightarrow} V_1
\stackrel{\varphi_1}{\dashrightarrow} \cdots
\stackrel{\varphi_{i-1}}{\dashrightarrow} V_i
\stackrel{\varphi_{i}}{\dashrightarrow} V_{i+1}
\stackrel{\varphi_{i+1}}{\dashrightarrow}
\cdots \stackrel{\varphi_{l-2}}{\dashrightarrow}
V_{l-1} \stackrel{\varphi_{l-1}}{\dashrightarrow} V_l = X_2
$$
de sorte que 
$\varphi = \varphi_{l-1} \circ \varphi_{l-2} \circ \cdots \varphi_1 \circ
\varphi_0$ et pour tout $i$,
$\varphi_i:V_{i} \dashrightarrow V_{i+1}$ ou $\varphi_i^{-1}:V_{i+1}
\dashrightarrow V_{i}$
est une application r\'eguli\`ere obtenue en \'eclatant une
sous-vari\'et\'e irr\'eductible lisse. 
\end{theo}

Ce texte est pour l'essentiel une reprise de 
\cite{AKMW99} et \cite{Mat99}, dans lesquels on trouvera 
beaucoup plus de pr\'ecisions, une liste de r\'ef\'erences
tr\`es compl\`ete ainsi qu'une discussion d\'etaill\'ee des
extensions ou g\'en\'eralisations du th\'eor\`eme \ref{faible general}.

Mentionnons aussi que le th\'eor\`eme de factorisation
admet les raffinements fondamentaux suivants :
{\em 
\begin{enumerate}
\item [$\bullet$] si $\varphi$ est un
isomorphisme sur un ouvert $U$, le centre de chaque $\varphi_i$
ou $\varphi_i^{-1}$ peut \^etre choisi disjoint de $U$,
\item [$\bullet$] si $X_1$ et $X_2$ sont des vari\'et\'es projectives,
chaque $V_i$ peut \^etre choisie projective,
\item [$\bullet$] si $X_1$ et $X_2$ ne sont pas suppos\'ees compl\`etes,
toute application birationnelle propre (voir
\cite{Har77} pour cette notion)
entre $X_1$ et $X_2$ se factorise en 
une suite d'\'eclatements et de contractions de centres lisses,
\item [$\bullet$] si $X_1$ et $X_2$ sont deux 
vari\'et\'es a\-nalytiques complexes compactes lisses, 
toute application bim\'ero\-mor\-phe propre 
entre $X_1$ et $X_2$ se factorise en 
une suite d'\'eclatements et de contractions de centres analytiques 
lisses.
\end{enumerate}
}

Pour simplifier l'exposition, nous nous contenterons de renvoyer
\`a \cite{AKMW99} et \cite{Mat99} pour une d\'e\-monstra\-tion 
de ces points.
 
Faisons aussi les remarques suivantes :
{\em 
\begin{enumerate}
\item [$\bullet$] la factorisation n'est \'evidemment pas unique,
\item [$\bullet$] dans le cas des surfaces, si $\varphi$
est r\'eguli\`ere, alors  
toutes les $\varphi_i$ peuvent \^etre choisies 
r\'eguli\`eres. Ceci est faux en dimension
sup\'erieure ou \'egale \`a trois, 
\item [$\bullet$] la factorisation d\'ecrite
dans le th\'eor\`eme \ref{faible general} est une 
solution positive au probl\`eme de factorisation faible~;
le probl\`eme de factorisation forte 
(au sens o\`u il existe un entier $i_0$
tel que pour tout $i\leq i_0$, $\varphi_{i} ^{-1}$
est r\'eguli\`ere et pour tout $i\geq i_0+1$, $\varphi_{i}$
est r\'eguli\`ere) est un probl\`eme ouvert \`a ce jour 
en dimension
sup\'erieure ou \'egale \`a trois.
\end{enumerate}
} 

\subsection{Quelques mots sur la d\'emonstration du th\'eor\`eme
\ref{faible general}}
La d\'emonstration du th\'eor\`eme
\ref{faible general} est en certain sens une victoire de la g\'eom\'etrie
torique. En effet, la ligne directrice est 
une r\'eduction (en plusieurs \'etapes, certaines utilisant
fondamentalement des techniques toriques) au cas d'une application
birationnelle \'equivariante entre vari\'et\'es toriques, 
cas r\'esolu par Morelli \cite{Mor96}
et W{\l}odarczyk \cite{Wlo97},
puis \'etendu au cadre toro\"{\i}dal par 
Abramovich, Matsuki et
Rashid \cite{AMR99}. Il est bien connu depuis quelques ann\'ees
que la g\'eom\'etrie des vari\'et\'es toriques
(ou des plongements toro\"{\i}daux), gouvern\'ee par 
des objets combinatoires simples issus de la g\'eom\'etrie convexe,
donne une bonne vision locale de certaines propri\'et\'es
des vari\'et\'es alg\'ebriques (voir par exemple l'existence 
d'alt\'erations
due \`a De Jong). L'expos\'e qui suit donnera, je l'esp\`ere, envie
au lecteur de mieux conna\^{\i}tre ou de d\'ecouvrir ces techniques ;
ce texte ne rempla\c cant certainement pas l'\'enorme 
effort p\'edagogique que constitue le texte \cite{Mat99}
de K.~Matsuki. 
   
\medskip

\subsection{Remerciements} 
Merci aux coll\`egues
qui m'\'ecoutent parler d'\'eclatements depuis plusieurs
ann\'ees, ils sont trop nombreux pour \^etre tous mentionn\'es. 
Merci \`a Michel Brion, Laurent Manivel, Kenji Matsuki et Emmanuel Peyre
pour m'avoir aid\'e \`a pr\'eparer ce texte, 
et plus particuli\`erement \`a St\'ephane Guillermou,
infatigable relecteur de multiples versions pr\'eliminaires.  
Je d\'edie ce travail \`a mon fils Alex, et \`a Marguerite sa maman.

\section{Cobordisme birationnel et action de $\KK ^*$}\label{cobord}

La notion de cobordisme birationnel a \'et\'e d\'egag\'ee par
W{\l}odarczyk \cite{Wlo97}, suite au travail fondamental de Morelli
\cite{Mor96} dans le cadre des vari\'et\'es toriques. L'id\'ee essentielle 
est la suivante : la th\'eorie de Morse sur les vari\'et\'es 
(diff\'erentiables r\'eelles) permet de reconstruire 
topologiquement une vari\'et\'e donn\'ee $X$ \`a partir d'une fonction de
Morse $f$ sur $X$. Le passage des points critiques de $f$
(qui sont aussi les points fixes du champ de
vecteurs
${\rm grad} (f)$ lorsque 
$X$ est munie d'une m\'etrique riemannienne)
correspond aux changements de topologie par ajout
d'une cellule. A
un morphisme birationnel entre $X_1$ et $X_2$ de dimension $n$, on 
associe une vari\'et\'e alg\'ebrique de dimension $n+1$
munie d'une action de $\KK ^*$, avec un ordre 
sur les composantes connexes de points fixes. 
L'action de $\KK^*$ joue le r\^ole du champ de vecteurs
${\rm grad} (f)$ et
\`a chaque
composante connexe de $\KK ^*$-points fixes, on associe une
application birationnelle \'el\'ementaire, dont on montre
qu'elle est ``localement torique''.  
  
\subsection{Rappels sur les $\KK^*$-actions}

Lorsque le groupe multiplicatif $\KK ^*$
agit alg\'ebriquement sur une vari\'et\'e alg\'ebrique $X$,
pour $t\in \KK ^*$ et $x \in X$, le r\'esultat de l'action
de $t$ sur $x$ sera not\'e $t \cdot x$.
On note $X/\KK ^*$ l'ensemble des orbites de
l'action et $X^{\KK^*}$ l'ensemble des $\KK^*$-points 
fixes de $X$.

\medskip

\noindent {\bf D\'efinition.}
Un {\em bon quotient} ou {\em quotient cat\'egorique}
$Y = X/\!/\KK^*$
est la donn\'ee d'une vari\'et\'e alg\'ebrique $Y$ et d'une application
r\'eguli\`ere $\pi : X \to Y$ constante sur les $\KK ^*$-orbites
de sorte que pour tout ouvert affine 
$U$ de $Y$, $\pi ^{-1}(U)$ est un ouvert affine de $X$
et l'application induite $\pi ^* : {\mathcal O}_Y (U)
\to ({\mathcal O}_X (\pi ^{-1}(U)))^{\KK ^*}$
est un isomorphisme (ici, $({\mathcal O}_X (\pi ^{-1}(U)))^{\KK ^*}$
d\'esigne les \'el\'ements $\KK ^*$-invariants 
de ${\mathcal O}_X (\pi ^{-1}(U))$).
Si de plus les fibres de $\pi$ sont exactement
les orbites, alors $Y$
est appel\'e {\em quotient g\'eom\'etrique}.

\medskip

Rappelons (voir par exemple \cite{MFK94} ou \cite{Dol94})
que $\KK ^*$ \'etant un groupe r\'eductif, 
pour toute vari\'et\'e affine $X$ munie d'une action 
alg\'ebrique de $\KK^*$, 
le quotient cat\'egorique
$X/\!/\KK^*$ existe et ses points correspondent aux $\KK ^*$-orbites
ferm\'ees. De plus $X/\!/\KK^*$ est normale si $X$ l'est.

\medskip

\noindent {\bf Notations.}
Soit $X$ une vari\'et\'e alg\'ebrique sur laquelle
$\KK ^*$
agit alg\'ebriquement.
Introduisons les deux sous-ensembles localement
ferm\'es de $X$ suivant : 
$$ X _+ = \{ x \in X \, | \, \lim_{t \to \infty} t\cdot x 
\,\, \mbox{n'existe pas dans}\, \, X \} $$ et
$$ X _- = \{ x \in X \, |\, \lim_{t \to 0} t\cdot x 
\, \,\mbox{n'existe pas dans}\, \, X \} .$$ 
Pr\'ecisons ici que ``$\lim_{t \to 0} t\cdot x$
existe dans $X$'' signifie que l'application r\'eguli\`ere 
de $\KK^*$ dans $X$ qui \`a $t$ associe $t\cdot x$
s'\'etend en une application r\'eguli\`ere de $\KK$
dans $X$ ayant $\lim_{t \to 0} t\cdot x$ pour valeur en $0$.

\medskip

\noindent {\bf Premier exemple fondamental.}
Consid\'erons l'action alg\'ebrique de $\KK ^*$
sur $X = \KK ^{n+1}$ d\'efinie par 
$t\cdot (x_1,\ldots,x_{n+1}) = (t^{a_1}x_1,\ldots,t^{a_{n+1}}x_{n+1})$
o\`u les $a_i$ sont des entiers premiers entre eux 
tels que $a_i <0$ pour $1\leq i \leq \alpha$
et $a_i > 0$ pour $\alpha +1 \leq i \leq n+1$
pour un entier $\alpha$ tel que $2 \leq \alpha \leq n$. 
Alors $X_+ = X \setminus (\KK ^{\alpha} \times \{0\})$
et $X_- = X \setminus (\{0\} \times \KK ^{n+1-\alpha})$.
Les quotients g\'eom\'etriques $X_+/\KK ^*$
et $X_-/\KK ^*$ existent, ainsi que le quotient cat\'egorique
$X/\!/\KK^*$ et on a un diagramme commutatif~:

\centerline{
\xymatrix{ X_-/\KK ^* \ar[rd]^{\varphi_-}\ar@{-->}[rr]^{\varphi} 
& & X_+/\KK ^* \ar[ld]_{\varphi_+}\\
 &X/\!/\KK^*& 
}
}

Si $\underline{0} \in X/\!/\KK^* $ d\'esigne l'unique
$\KK^*$-orbite ferm\'ee de $X$, la fibre $\varphi_+ ^{-1}(\underline{0})$ 
(resp. $\varphi_- ^{-1}(\underline{0})$)
est isomorphe \`a l'espace projectif \`a poids
$\PP (a_{\alpha +1},\ldots,a_{n+1})$ (resp.    
$\PP (a_{1},\ldots,a_{\alpha})$)
De plus, $\varphi_+$ et $\varphi_-$ se restreignent
en des isomorphismes sur l'ouvert
$(X_- \cap X_+)/\KK ^*$, ainsi $\varphi = \varphi_+^{-1}\circ 
\varphi_- : X_-/\KK ^* \dra X_+/\KK ^*$ est une application
birationnelle.

\subsection{Cobordisme birationnel}

La d\'efinition suivante sera l'outil fondamental 
dans toute la suite.

\medskip 

\noindent {\bf D\'efinition.} 
Soit $\varphi:X_1 \dashrightarrow X_2$ 
une application birationnelle entre deux vari\'et\'es 
alg\'ebriques compl\`etes et lisses
$X_1$ et $X_2$ sur $\KK$.
Un {\em cobordisme birationnel} pour $\varphi$ 
est la donn\'ee d'une vari\'et\'e alg\'ebrique normale
$B$ telle que :
\begin{enumerate}
\item [(i)] le groupe multiplicatif $\KK ^*$
agit de fa\c con effective sur $B$ ({\em i.e.}
$\bigcap_{x \in B} \Stab (x) = 1$),
\item [(ii)] les ensembles $B_-$ et $B_+$
sont des ouverts (de Zariski) non vides,
\item [(iii)] les quotients g\'eom\'etriques 
$B_-/\KK ^*$ et $B_+/\KK ^*$ existent
et sont respectivement isomorphes \`a $X_1$
et $X_2$, 
\item [(iv)] si $\psi: B_- \dra B_+$ est l'application birationnelle
induite par les inclusions $B_-\cap B_+ \subset B_-$ et
$B_-\cap B_+ \subset B_+$, le diagramme suivant est commutatif :

\centerline{
\xymatrix{ B_- \ar[d] \ar@{-->}[r]^{\psi}& B_+ \ar[d] \\
X_1 \ar@{-->}[r]^{\varphi} &X_2
}
}
\end{enumerate}

\medskip

Dans le premier exemple fondamental,
$X=\KK ^{n+1} $ est un cobordisme 
birationnel pour $\varphi : X_-/\KK ^* \dra X_+/\KK ^*$.  

\medskip

\noindent {\bf Deuxi\`eme exemple fondamental (voir aussi
\cite{Ful84}).}
Soient $X$ une vari\'et\'e alg\'ebrique compl\`ete et lisse, 
$Y$ une sous-vari\'et\'e irr\'eductible lisse de $X$
et $\varphi : X_1 \to X_2=X$ l'\'eclatement de $X$
de centre $Y$.
Construisons un cobordisme birationnel pour
$\varphi$.  
Soit $W$ la vari\'et\'e alg\'ebrique
$X_2 \times \PP ^1$ sur laquelle
le groupe $\KK ^*$ agit alg\'ebriquement par
multiplication sur le second facteur.
Notons 
$\overline{Y} = Y \times \{0\} \subset W$
et soit $\overline{B}$ la vari\'et\'e alg\'ebrique compl\`ete et lisse
obtenue en \'eclatant $W$ le long de $\overline{Y}$. 
Comme $\overline{Y}$ est inclus dans l'ensemble des
$\KK ^*$-points fixes de $W$, l'action de $\KK ^*$
sur $W$ se rel\`eve en une action sur $\overline{B}$.
La transform\'ee stricte $D_1$ de $X_2 \times \{ 0 \}$ est isomorphe
\`a $X_1$ si bien que $\overline{B}$ poss\`ede deux diviseurs 
constitu\'es de $\KK ^*$-points fixes, l'un, $D_1$, isomorphe \`a $X_1$ et 
l'autre, not\'e $D_2$ et \'egal \`a l'image inverse dans $\overline{B}$ de
$X_2 \times \{\infty \}$, isomorphe \`a $X_2$.
Posons alors $B = \overline{B} \setminus (D_1 \cup D_2)$.
Soit $E$ le diviseur exceptionnel de l'\'eclatement
$\overline{B} \to W$. Ce diviseur 
est isomorphe \`a $\PP (N_{\overline{Y}/W})
= \PP( N_{Y/X_2} \oplus \underline{0})$
et est $\KK ^*$-invariant
($\underline{0}$ d\'e\-si\-gne le fibr\'e en droites trivial).
Remarquons que l'ensemble des $\KK ^*$-points fixes de $B$ correspond
\`a l'image de $\PP (N_{X_2 \times \{0\}/W |\overline{Y}})$
dans l'identification pr\'ec\'edente
et est donc naturellement isomorphe \`a $Y$.
C'est ensuite un exercice facile de voir 
que $B$ est un cobordisme birationnel pour $\varphi$.
La figure suivante repr\'esente le cas de l'\'eclatement 
d'un point dans une surface et devrait \'eclairer le lecteur.

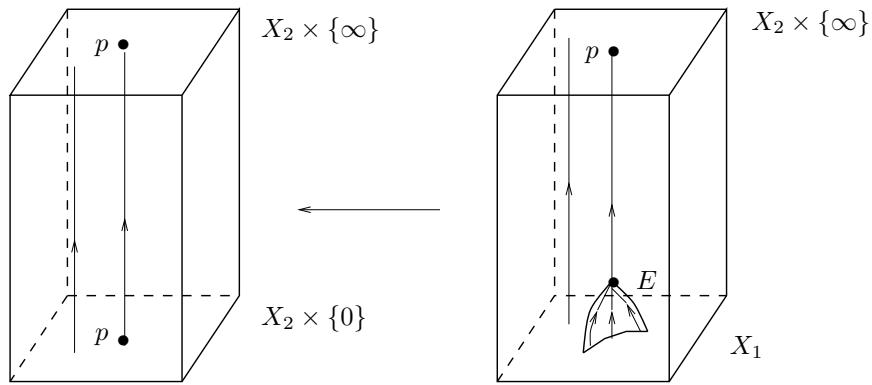
\begin{figure}[hbtp]
  \begin{center}
    \leavevmode
    \input{fig1.pstex_t}
    \caption{Cobordisme birationnel : \'eclatement d'un point}
    \label{}
  \end{center}
\end{figure}

\subsection{Construction de cobordisme birationnel}\label{existe cobord} 

On d\'emontre ici le r\'esultat suivant, d\^u \`a W{\l}odarczyk 
\cite{Wlo99} :

\begin{theo}\label{existence cobordisme}
Soit $\varphi : X_1 \to X_2$ une application r\'eguli\`ere
birationnelle entre deux vari\'et\'es projectives lisses. 
Alors, il existe une vari\'et\'e projective lisse
$\overline{B}$ munie d'une action alg\'ebrique effective
de $\KK ^*$ v\'erifiant :
\begin{enumerate}
\item [(i)] il existe deux plongements 
$\iota _1 : X_1 \to \overline{B}^{\KK^*}$
et $\iota _2 : X_2 \to \overline{B}^{\KK^*}$
d'images disjointes,
\item [(ii)] la vari\'et\'e $B = \overline{B} \setminus (\iota_1(X_1)
\cup \iota_2(X_2))$ est un cobordisme birationnel pour $\varphi$.
\end{enumerate}
\end{theo}

\noindent {\em D\'emonstration.}
Elle suit de pr\`es la construction 
du deuxi\`eme exemple fondamental.
Comme $\varphi : X_1 \to X_2$ est une application r\'eguli\`ere
birationnelle entre deux vari\'et\'es projectives, il 
existe un faisceau d'id\'eaux ${\mathcal I} \subset {\mathcal O}_{X_2}$
coh\'erent tel que $\varphi$ soit l'\'eclatement de 
${\mathcal I}$ (voir par exemple \cite{Har77} page 166).
Soit $W$ la vari\'et\'e alg\'ebrique
$X_2 \times \PP ^1$ sur laquelle
le groupe $\KK ^*$ agit alg\'ebriquement par
multiplication sur le second facteur.
Soit ${\mathcal J}= 
(p_1^{-1}{\mathcal I}+p_2^{-1}{\mathcal I}_0){\mathcal O}_W$
o\`u ${\mathcal I}_0$ est l'id\'eal du point $0$ de $\PP^1$ et 
soit $\pi : \tilde W \to W$
l'\'eclatement de ${\mathcal J}$. Comme ${\mathcal J}$ 
est $\KK^*$-invariant, l'action de $\KK^*$
sur $W$ se rel\`eve \`a $\tilde W$.  
La vari\'et\'e $\tilde W$ est projective et en g\'en\'eral singuli\`ere.
La transform\'ee stricte $D_1$ de $X_2 \times \{ 0 \}$ est isomorphe
\`a $X_1$ et, \'etant lisse, de Cartier et d'\'equation
locale d'ordre $1$, est incluse dans le lieu 
r\'egulier de $\tilde W$. 
Soit alors $\overline{B}$ une d\'esingularisation
canonique 
de $\tilde W$ (une telle d\'esingularisation
est obtenue par une suite d'\'eclatements le long de 
centres lisses disjoints du lieu r\'egulier de $\tilde W$
et est naturellement munie d'une action alg\'ebrique de $\KK^*$
relevant l'action de $\KK^*$ sur $W$ ; de telles d\'esingularisations
existent, voir la partie \ref{desing} pour plus de d\'etails).
Alors, la pr\'e-image $D'_1$ de $D_1$ dans $\overline{B}$
est isomorphe \`a $X_1$, la pr\'e-image $D'_2$ de
$X_2 \times \{\infty\}$ dans $\overline{B}$
est isomorphe \`a $X_2$ et $B = \overline{B} \setminus (D'_1\cup D'_2)$
est un cobordisme birationnel pour $\varphi$.
\finpreuve

\subsection{Filtrabilit\'e}\label{filtre}
La notion de cobordisme filtrable est due \`a Morelli
\cite{Mor96} et W{\l}odarczyk \cite{Wlo97}, 
son origine se situe
dans les travaux de Bialynicki-Birula.

\medskip

\noindent {\bf D\'efinition.} 
Soit $B$ une vari\'et\'e alg\'ebrique munie d'une action 
alg\'ebrique de $\KK ^*$.
Si $F_1$ et $F_2$ sont deux composantes connexes
de $B^{\KK^*}$, on dit que 
{\em $F_1$ pr\'ec\`ede $F_2$}, que l'on note 
$F_1 \prec F_2$, s'il existe $x\in B \setminus B^{\KK^*}$
tel que $\lim_{t \to 0} t\cdot x \in F_1$ et
$\lim_{t \to \infty} t\cdot x \in F_2$.

\medskip

La relation $\prec$ n'est en g\'en\'eral pas transitive. 

\medskip

\noindent {\bf D\'efinition.} 
Soit $B$ une vari\'et\'e alg\'ebrique munie d'une action 
alg\'ebrique de $\KK ^*$. On dit que $B$ est {\em filtrable}
si'il
n'y a pas de $\prec$-cycle de composantes connexes de $B^{\KK^*}$
$$ F_1 \prec F_2 \prec \cdots \prec F_m \prec F_1.$$
En particulier, il n'y a pas de composante connexe $F$ de
$B^{\KK^*}$ telle que $F \prec F$.

\medskip

Le lemme suivant est \'el\'ementaire mais essentiel :

\begin{lemm}\label{lemme filtre}
Soit $B$ une vari\'et\'e alg\'ebrique lisse munie d'une action 
alg\'ebrique de $\KK ^*$. Si $B$ est quasi-projective,
alors $B$ est filtrable. 
\end{lemm}

\noindent {\em D\'emonstration.}
Par un r\'esultat de Sumihiro \cite{Sum74} \cite{Sum75}, 
il existe une 
immersion localement ferm\'ee et \'equivariante 
de $B$ dans un espace projectif $\PP (V)$
o\`u $V$ est un espace vectoriel sur lequel $\KK^*$ agit
lin\'eairement et alg\'ebriquement. Comme l'action de $\KK ^*$
est alg\'ebrique, $V$ se d\'ecompose en espaces propres
$V =\bigoplus _{k=0}^{m} V(a_k)$ o\`u
les $a_k$ sont des entiers relatifs ordonn\'es $a_0 < \cdots < a_m$
et 
$\KK^*$ agit sur $x_k \in V(a_k)$ par $t\cdot x_k = t^{a_k}x_k$
(les $a_k$ sont les {\em poids} de la repr\'esentation $V$).
L'observation suivante est \'el\'ementaire mais cruciale :
{\em soit $x \in \PP(V)$ 
que l'on rel\`eve en $\bar x = x_0 \oplus \cdots \oplus x_m$
dans $V$. Soit ${\rm min} (x)$ (resp. ${\rm max} (x)$)
le plus petit (resp. grand) indice $i$ dans $\{ 0,\ldots,m\}$
tel que $x_i\neq 0$. Alors $\lim _{t \to 0} t \cdot x$
(resp. $\lim _{t \to \infty} t \cdot x$)
est l'image dans $\PP (V)$ de $x_{{\rm min} (x)}$ (resp. 
$x_{{\rm max} (x)}$). 
}

\noindent Les composantes connexes de $(\PP (V))^{\KK ^*}$
sont les $C_{a_k}= \PP (V(a_k))$ et le lemme repose sur 
le fait suivant : 
s'il existe $x$ dans $\PP (V) \setminus (\PP (V))^{\KK ^*}$
tel que  
$\lim_{t \to 0} t\cdot x \in C_{a_i}$ et
$\lim_{t \to \infty} t\cdot x \in C_{a_j}$, alors, d'apr\`es l'observation
pr\'ec\'edente,
$a_i < a_j$.
De l\`a, si $F$ est une composante connexe de
$B^{\KK^*}$ et si $a(F)$ est l'unique entier
tel que $F \subset C_{a(F)}$, on d\'eduit que
$F \prec F'$ implique $a(F) < a(F')$.
Il s'ensuit que $B$ est filtrable. \finpreuve  

\subsection{D\'ecomposition d'un cobordisme birationnel et 
th\'eorie g\'eom\'etrique
des invariants}\label{GIT}
Soit $\varphi : X_1 \to X_2$ une application r\'eguli\`ere
birationnelle entre deux vari\'et\'es projectives lisses
et soit $B$ un cobordisme birationnel quasi-projectif 
(et donc filtrable d'apr\`es ce qui pr\'ec\`ede) pour $\varphi$ 
donn\'e par le th\'eor\`eme \ref{existence cobordisme}.
Dans le paragraphe pr\'ec\'edent, on a associ\'e \`a 
$B$ une suite de poids entiers ordonn\'es $(a_i)_{i=0,\ldots,m}$,
et quitte \`a remplacer $V$ par ${\rm Sym}^2 (V)$,
on peut supposer de plus que tous les poids $a_k$
sont pairs, en particulier
$a_i<a_i +1 < a_{i+1}$.
Rappelons aussi qu'\`a  
toute
composante connexe $F$ de $B^{\KK^*}$
correspond l'un de ces poids  
$a(F)$. Si $x$ appartient \`a $B^{\KK^*}$,
$x$ appartient \`a une unique composante connexe
$F$ de $B^{\KK^*}$ et   
on note $a(x)$ le poids $a(F)$ correspondant.

\medskip

\noindent {\bf Notations.}  
Soit $a_i$ l'un des poids pr\'ec\'edents. 
On note $B_{a_i}$ le compl\'ementaire dans $B$ de
$$\{ x \, | x_{\infty} = \lim_{t \to \infty} t\cdot x 
\mbox{ existe dans } B  \mbox{ et } 
a(x_{\infty}) < a_i \} 
\, \, \cup $$
$$\{ x \, | x_0= \lim_{t \to 0} t\cdot x 
 \mbox{ existe dans }  B  \mbox{ et } 
a(x_0) > a_i \}.
$$

Remarquons que chaque $B_{a_i}$ contient l'ouvert $B_- \cap B_+$
et que si $x$ est un point fixe dans $B_{a_i}$, alors 
$a(x) = a_i$. En particulier, d'apr\`es le paragraphe \ref{filtre},
il n'existe pas de point $x$ dans $B_{a_i}\setminus B_{a_i}^{\KK ^*}$
tel que $\lim_{t \to 0} t \cdot x$ et 
$\lim_{t \to \infty} t \cdot x$ existent dans $B_{a_i}$.

Le lemme suivant est imm\'ediat~:

\begin{lemm}
Avec les notations 
pr\'ec\'edentes, 
$(B_{a_i})_+ = (B_{a_{i+1}})_-$
pour tout $0 \leq i \leq m-1$. De plus,
$B_- = (B_{a_0})_-$ et $B_+ = (B_{a_m})_+$.
\end{lemm}

Le r\'esultat suivant est d\^u \`a W{\l}odarczyk \cite{Wlo99},
il voit son origine dans les travaux de 
Guil\-lemin-Sternberg \cite{GuS89}
et Brion-Procesi \cite{BrP90} :

\begin{prop}Le quotient 
cat\'egorique $B_{a_i} /\!/\KK^*$ et les quotients
g\'eom\'etriques $(B_{a_i})_-/\KK^*$ et
$(B_{a_i})_+/\KK^*$ existent.
De plus, $B_{a_i}$ est un cobordisme
birationnel pour $$\varphi _i:(B_{a_i})_-/\KK^* \dra (B_{a_i})_+/\KK^*.$$
\end{prop}

\noindent {\em D\'emonstration.}
Les notations \'etant celles de la d\'emonstration
du lemme \ref{lemme filtre}, soit $\rho _0$
l'action de $\KK ^*$ sur 
$V =\bigoplus _{k=0}^{m} V(a_k)$ et pour un entier
$r \in \ZZ$, soit $\rho_r$
l'action de $\KK ^*$ obtenue en ``tordant''
$\rho _0$ par $t^{-r}$. Autrement dit, si
$x_k \in V(a_k)$, on a $\rho_r(t) \cdot x_k = t^{a_k-r}x_k$.
Evidemment, l'action de $\rho_r$ 
sur $\PP (V)$ est \'egale \`a l'action initiale $\rho_0$,
mais $\rho_r$ induit un changement de lin\'earisation 
sur le $\KK ^*$-fibr\'e en droites ample ${\mathcal O}_{\PP (V)}(1)$.
Sous l'hypoth\`ese que les poids de l'action v\'erifient
$a_i +1 < a_{i+1}$, on observe que $B_{a_i}$ (resp. $(B_{a_i})_-$,
resp. $(B_{a_i})_+$)  
est le lieu des points semi-stables de $B$
pour le fibr\'e en droites ${\mathcal O}_{\PP (V)}(1)$ lin\'earis\'e par 
$\rho_{a_i}$ (resp. $\rho_{a_i -1}$, resp. $\rho_{a_i +1}$).  
La th\'eorie g\'eom\'etrique des invariants \cite{MFK94} \cite{Dol94}
assure alors que 
les quotients cat\'egoriques $B_{a_i} /\!/\KK^*$, $(B_{a_i})_-/\!/\KK^*$
et $(B_{a_i})_+/\!/\KK^*$ existent, et comme de plus 
chaque orbite de $(B_{a_i})_-$
(resp. $(B_{a_i})_+$) est ferm\'ee dans $(B_{a_i})_-$ (resp. $(B_{a_i})_+$),
les quotients $(B_{a_i})_-/\!/\KK^*$
et $(B_{a_i})_+/\!/\KK^*$ sont des quotients g\'eom\'etriques.
De plus, le quotient cat\'egorique 
$\pi_i : B_{a_i} \to B_{a_i} /\!/\KK^*$ est un morphisme affine.
Notons enfin qu'il y a un diagramme commutatif~:
 
\centerline{
\xymatrix{ (B_{a_i})_-/\KK^* \ar[rd]_{(\varphi _i)_-}  
\ar@{-->}[rr]^{\varphi_i}&
& (B_{a_i})_+/\KK^*  \ar[ld]^{(\varphi _i)_+} \\
& B_{a_i}/\!/\KK^* & 
}
}   
\noindent si bien que le cobordisme birationnel 
$\varphi _i:(B_{a_i})_-/\KK^* \dra (B_{a_i})_+/\KK^*$ s'interpr\`ete
comme un changement de lin\'earisation pour
la restriction \`a $B_{a_i}$ du
$\KK ^*$-fibr\'e en droites ample ${\mathcal O}_{\PP (V)}(1)$. 
\finpreuve

\bigskip

\noindent 
{\bf Faisons le point :} si $\varphi : X_1 \to X_2$ est une application 
r\'eguli\`ere
birationnelle entre deux vari\'et\'es projectives lisses,
ce qui pr\'ec\`ede montre que 
$\varphi$ se factorise en 
une suite d'applications birationnelles 
$$ X_1 = V_0 \stackrel{\varphi_0}{\dashrightarrow} V_1
\stackrel{\varphi_1}{\dashrightarrow} \cdots
\stackrel{\varphi_{i-1}}{\dashrightarrow} V_i
\stackrel{\varphi_{i}}{\dashrightarrow} V_{i+1}
\stackrel{\varphi_{i+1}}{\dashrightarrow}
\cdots \stackrel{\varphi_{m-1}}{\dashrightarrow}
V_{m} \stackrel{\varphi_{m}}{\dashrightarrow} V_{m+1} = X_2
$$
o\`u chaque 
$\varphi_i:V_{i} \dashrightarrow V_{i+1}$ 
est un cobordisme
birationnel 
$\varphi _i:(B_{a_i})_-/\KK^* \dra (B_{a_i})_+/\KK^*$,
la vari\'et\'e $B_{a_i}$ \'etant une $\KK ^*$-vari\'et\'e quasi-projective
pour laquelle le quotient ca\-t\'e\-go\-rique $B_{a_i}/\!/\KK^*$
existe.

\subsection{Structures localement toriques}\label{slt}
Nous allons expliquer dans ce paragraphe, en montrant
que l'on peut munir les
vari\'et\'es $B_{a_i}$ construites pr\'ec\'edemment
d'un atlas
de ``cartes \'etales fortement toriques'',
que la d\'e\-com\-po\-sition
obtenue pr\'ec\'edemment est une d\'e\-com\-po\-si\-tion
en applications ``localement toriques''. 
Ceci est un point essentiel
pour la suite o\`u toutes les constructions \'elabor\'ees (en
particulier dans
la partie \ref{struc toroid}) et la v\'erification de
leurs propri\'et\'es se feront \`a l'aide de cet atlas.

\medskip

Rappelons qu'une application r\'eguli\`ere $f: Z \to X$ entre vari\'et\'es 
alg\'ebriques est {\em \'etale} si elle est lisse de dimension relative $0$.
De fa\c con \'equivalente, $f$ est \'etale si pour
tout $z \in Z$, $f$ induit
un isomorphisme $f^* : \hat{{\mathcal O}}_{X,f(z)} \to
\hat{{\mathcal O}}_{Z,z}$ entre les compl\'et\'es des anneaux locaux. 
Sur le corps des nombres complexes, $f$
est \'etale si et seulement si $f$ est un biholomorphisme local.

\medskip

\noindent {\bf D\'efinition.} Soit 
$f : Z \to X$ un morphisme \'etale, $\KK ^*$-\'equivariant
entre deux vari\'et\'es alg\'ebriques affines 
munies d'une action de $\KK ^*$.
On dit que $f$ est {\em fortement \'etale} si 
et seulement si $f/\!/\KK ^* : Z/\!/\KK ^* \to X/\!/\KK ^*$
est \'etale et l'application naturelle
$Z \to Z/\!/\KK ^* \times_{X/\!/\KK ^*} X$ est un isomorphisme.

\medskip

Soit $V$ une vari\'et\'e alg\'ebrique lisse munie d'une action 
alg\'ebrique de $\KK ^*$ et
soit $x\in V$ un $\KK^*$-point fixe.
Par un r\'esultat de Sumihiro \cite{Sum74} \cite{Sum75},
il existe un voisinage $\KK^*$-invariant et affine $W_{x}$ de $x$ 
dans $V$.
Comme $x$ est un point fixe, le lemme
fondamental de Luna \cite{Lun73} \cite{MFK94} 
assure qu'il existe un voisinage 
$\KK^*$-invariant et affine $V_x \subset W_{x}$ de $x$,
satur\'e pour la projection $\pi_x : W_x \to W_x/\!/\KK^*$,
et un morphisme $\KK ^*$-\'equivariant et fortement \'etale 
$\eta_x : V_x \to X_x := T_x V$ (o\`u
$T_x V$ est l'espace tangent de $V$ en $x$).
Comme l'action de $\KK^*$ sur 
$T_{x}V$ se diagonalise, il existe 
une base de vecteurs propres
de $T_{x}V$, faisant
de $T_{x}V$ une vari\'et\'e torique affine lisse 
o\`u $\KK ^*$ agit comme sous-groupe \`a un param\`etre du tore
(nous renvoyons au paragraphe suivant pour
les notions de g\'eom\'etrie torique).

\medskip

En adaptant un peu cette construction, on obtient la
proposition suivante, o\`u l'on a repris les notations
des paragraphes pr\'ec\'edents :

\begin{prop}\label{cartes toriques}
Soit $B_{a_i}$ comme dans le paragraphe \ref{GIT} et soit
$\pi_i : B_{a_i}\to B_{a_i}/\!/\KK^*$
le quotient cat\'ego\-ri\-que.
Alors, pour tout 
$x \in B_{a_i}$, 
il y a un voisinage  
$\KK^*$-invariant et affine $V_x$ de $x$, satur\'e pour $\pi_i$,
une vari\'et\'e torique affine lisse $X_x$ 
et un morphisme $\KK ^*$-\'equivariant et fortement \'etale 
$\eta_x : V_x \to X_x$.
\end{prop}

Mentionnons ici que l'on peut choisir les $V_x$ satur\'es pour 
$\pi_i$ car $\pi_i$ est un morphisme affine.
Mentionnons aussi que la version
en caract\'eristique non nulle du r\'esultat de Luna 
donn\'ee par Bardsley et Richardson \cite{BaR85} peut s'appliquer ici
si on ne suppose plus $\KK$ de caract\'eristique nulle.

\medskip

\noindent {\bf D\'efinition.} Dans la proposition \ref{cartes toriques},
pour $x \in B_{a_i}$, on appelle {\em carte torique fortement \'etale
en $x$} la donn\'ee de $\eta_x : V_x \to X_x$.

\medskip

Ces cartes toriques fortement \'etales ont des propri\'et\'es
tr\`es fortes, en particulier :
\begin{enumerate}
\item [$\bullet$] $(B_{a_i})_- \cap V_x = (V_x)_-$
et $(B_{a_i})_+ \cap V_x = (V_x)_+$,
\item [$\bullet$] le morphisme $\eta_x$ se restreint en des morphismes 
$(\eta_x)_- :=(\eta_x) _{|(V_x)_-} : 
(V_x)_- \to (X_x)_-$ et 
$(\eta_x)_+ :=(\eta_x) _{|(V_x)_+} :
(V_x)_+ \to (X_x)_+$ 
fortement \'etales et il y a un diagramme commutatif~:  

\centerline{
\xymatrix{ (V_x)_-/\KK^*  \ar[d]_{(\eta_x)_-/\KK^*}   \ar[rd] & 
& (V_x)_+/\KK^*  \ar[d]^{(\eta_x)_+/\KK^*} \ar[ld] \\
(X_x)_-/\KK^*  \ar[rd]  & V_x/\!/\KK^* \ar[d] 
& (X_x)_+/\KK^* \ar[ld] \\
&X_x/\!/\KK^* & 
}
}
\end{enumerate}

\medskip

La proposition \ref{cartes toriques} 
signifie donc que le cobordisme
birationnel $\varphi _i:(B_{a_i})_-/\KK^* \dra (B_{a_i})_+/\KK^*$
est {\em localement torique} au sens o\`u, au voisinage
de tout point, et \`a morphisme \'etale pr\`es, l'application
birationnelle $\varphi _i$ est torique.
Il est essentiel de remarquer que nous sommes encore loin
d'avoir d\'emontr\'e le th\'eor\`eme de factorisation :
d'une part les $(B_{a_i})_-/\KK^*$ et $(B_{a_i})_+/\KK^*$
sont des vari\'et\'es singuli\`eres en g\'en\'eral, et d'autre
part la structure torique d\'epend du point $x$ choisi
dans $B_{a_i}$. 

Le but de la partie \ref{struc toroid} est 
de d\'ecomposer un cobordisme birationnel donn\'e
en une suite de cobordismes birationnels (en g\'en\'eral
singuliers) 
``toro\"{\i}daux'', ce qui signifie en un certain sens 
que la ``structure torique'' ne d\'ependra plus du point $x$
choisi.  

En attendant, dans la partie \ref{morelli} suivante,
nous nous int\'eressons au cadre torique. 

\section{Le th\'eor\`eme de Morelli et W{\l}odarczyk}\label{morelli}

Nous montrons ici le th\'eor\`eme de
factorisation faible des applications birationnelles
toriques.

\subsection{Rappels de g\'eom\'etrie torique et \'enonc\'e du th\'eor\`eme}
Des r\'ef\'erences usuelles de g\'eom\'etrie torique sont
\cite{Ewa96}, \cite{Ful93} et \cite{Oda88}.

\medskip

\noindent {\bf D\'efinition.} Une {\em vari\'et\'e torique} (de dimension
$n$) 
est une vari\'et\'e alg\'ebrique normale, contenant le
tore $T := (\KK ^*)^n$
comme ouvert (de Zariski), munie d'une action alg\'ebrique
de $T$ prolongeant l'action de $T$ sur lui-m\^eme.  

\medskip

Si $M$ est un r\'eseau 
({\em i.e.} un groupe ab\'elien 
libre) de rang $n$, une vari\'et\'e torique $X_{\Sigma}$
de dimension $n$ munie d'une action du tore $T = \Ho (M, \KK ^*)$
est d\'efinie par la donn\'ee d'un \'eventail
$\Sigma$, subdivision de l'espace vectoriel dual
$N_{\QQ}:= \Ho (M,\ZZ)\otimes_{\ZZ} \QQ$ par des c\^ones rationnels
poly\'edraux. 

Si $X_{\Sigma}$ est une vari\'et\'e torique d'\'eventail $\Sigma$,
\`a chaque c\^one $\sigma$ de $\Sigma$ correspond naturellement 
un ouvert
affine $T$-invariant $U_{\sigma}$ de $X_{\Sigma}$ de sorte que 
$X_{\Sigma}=\bigcup _{\sigma \in \Sigma} U_{\sigma}$.
De plus, si $\sigma \subset \tau$ ({\em i.e.} $\sigma$
est une face de $\tau$), alors $U_{\sigma} \subset U_{\tau}$.
Rappelons aussi qu'il y a une correspondance bijective entre les 
orbites de $T$ de codimension $r$ dans $X_{\Sigma}$ 
et les c\^ones de dimension $r$ de $\Sigma$. Pour $\sigma \in \Sigma$,
on note $V(\sigma)$ l'adh\'erence de l'orbite correspondant \`a $\sigma$ ;
si $\sigma \subset \tau$, alors $V(\tau) \subset V(\sigma)$.
Rappelons enfin 
qu'une 
vari\'et\'e torique est lisse si et seulement si
tous 
les c\^ones de $\Sigma$ sont {\em non-singuliers} ({\em i.e.} 
engendr\'es
par une famille de $N:=\Ho (M,\ZZ)$ 
pouvant se compl\'eter en une base de $N$)
et compl\`ete
si et seulement si le support de $\Sigma$ est $ N_{\QQ}$.

Le th\'eor\`eme de factorisation faible dans le cadre
torique est d\^u \`a Morelli \cite{Mor96} et W{\l}o\-darczyk
\cite{Wlo97}, voir
aussi \cite{AMR99} et \cite{Mat00} :

\begin{theo}\label{factorisation torique}
Soit $\varphi:X_1 \dashrightarrow X_2$ 
une application birationnelle \'equivariante entre deux vari\'et\'es 
toriques compl\`etes et lisses
$X_1$ et $X_2$ sur $\KK$. 
Alors, $\varphi$ se factorise en 
une suite d'\'eclatements et de contractions de centres lisses
invariants.
Autrement dit, il y a une suite 
d'applications birationnelles \'equivariantes entre vari\'et\'es toriques
compl\`etes et lisses
$$ X_1 = V_0 \stackrel{\varphi_0}{\dashrightarrow} V_1
\stackrel{\varphi_1}{\dashrightarrow} \cdots
\stackrel{\varphi_{i-1}}{\dashrightarrow} V_i
\stackrel{\varphi_{i}}{\dashrightarrow} V_{i+1}
\stackrel{\varphi_{i+1}}{\dashrightarrow}
\cdots \stackrel{\varphi_{l-2}}{\dashrightarrow}
V_{l-1} \stackrel{\varphi_{l-1}}{\dashrightarrow} V_l = X_2
$$
de sorte que 
$\varphi = \varphi_{l-1} \circ \varphi_{l-2} \circ \cdots \varphi_1 \circ
\varphi_0$ et pour tout $i$,
$\varphi_i:V_{i} \dashrightarrow V_{i+1}$ ou $\varphi_i^{-1}:V_{i+1}
\dashrightarrow V_{i}$
est une application r\'eguli\`ere obtenue en \'eclatant une
adh\'erence d'orbite du tore. 
\end{theo}

Dans la suite, nous donnons les grandes lignes de
la preuve du th\'eor\`eme \ref{factorisation torique}.
Mentionnons ici que malgr\'e les travaux \cite{Mor96}
et \cite{AMR99}, le probl\`eme de factorisation forte torique
est toujours ouvert, m\^eme en dimension trois (voir \cite{Mat00}
pour une discussion des lacunes de \cite{Mor96}
et \cite{AMR99}).

\subsection{Cobordisme torique en termes d'\'eventails, d'apr\`es Morelli}

Soient $M$ un r\'eseau de rang $n$ et $N$ le r\'eseau dual $\Ho (M,\ZZ)$. 
Soient $N^+$ le
r\'eseau $(n+1)$-dimensionnel $N \oplus \ZZ$,
$N_{\QQ} := N \otimes _{\ZZ} \QQ$, $N_{\QQ}^+ := N^+ \otimes _{\ZZ} \QQ$
et $\pi$ la projection $N_{\QQ}^+ \to N_{\QQ}$. 

Dans toute la suite 
et sauf mention explicite du contraire,
tous les c\^ones rationnels $\sigma$ de $N^+_{\QQ}$ que nous 
consid\`ererons 
seront suppos\'es {\em simpliciaux} ({\em i.e.}
engendr\'es par une famille libre d'\'el\'ements de $N^+$)
et {\em $\pi$-strictement convexes} ({\em i.e.} $\pi (\sigma)$ est un
c\^one strictement convexe de $N_{\QQ}$).

\medskip

\noindent {\bf D\'efinition.}
Soit $\sigma$ un c\^one rationnel de $N^+_{\QQ}$
On dit que $\sigma$ est {\em $\pi$-ind\'ependant} si
$\pi _{| \sigma}$ est injective (ceci signifie
que ${\rm Vect} (\sigma)$ ne contient pas la direction verticale
$\{ 0\} \oplus \QQ$). On dit que $\sigma$ est {\em $\pi$-d\'ependant} s'il
n'est pas $\pi$-ind\'ependant.
Soit $\sigma$ un c\^one rationnel $\pi$-ind\'ependant
de $N^+_{\QQ}$. On dit que $\sigma$ est {\em $\pi$-non-singulier}
si $\pi (\sigma)$ est un
c\^one non-singulier de $N_{\QQ}$. 
On dit que $\sigma$ est {\em $\pi$-singulier}
s'il
n'est pas $\pi$-non-singulier.
On dit qu'un \'eventail de $N^+_{\QQ}$ est {\em $\pi$-non-singulier}
si tous ses c\^ones $\pi$-ind\'ependants sont $\pi$-non-singuliers.

\medskip

Nous noterons 
$\nu = (0,1) \in N \oplus \ZZ = N^+$, ce vecteur correspond
\`a un sous-groupe \`a un param\`etre $\lambda _{\nu}$ 
du tore $N^+ \otimes _{\ZZ} \KK ^* = T \times \KK^*$. 

Soit $\Sigma$ un \'eventail simplicial de $N_{\QQ}^+$. D\'efinissons
les {\em faces sup\'erieure et inf\'erieure de $\Sigma$}
$$ \partial _+ (\Sigma) = \{ x \in \Sigma \, | \, x + \varepsilon \nu \notin 
\Sigma \mbox{ pour } \varepsilon > 0 \mbox{ petit} \}$$
et 
$$ \partial _- (\Sigma) = \{ x \in \Sigma \, | \, x - \varepsilon \nu \notin 
\Sigma \mbox{ pour } \varepsilon > 0 \mbox{ petit} \}.$$

\medskip

\noindent {\bf D\'efinition.}
Soit $\sigma$ un c\^one rationnel simplicial de $N^+_{\QQ}$.
Nous dirons que $\sigma$ est un {\em circuit} si
$\sigma$ est $\pi$-d\'ependant et
si toutes les faces strictes de 
$\sigma$ sont $\pi$-ind\'ependantes.

\medskip

Remarquons que tout c\^one rationnel
simplicial $\pi$-d\'ependant
de $N^+_{\QQ}$ contient ({\em i.e.} comme
face) un 
unique circuit.

Soit $\tau$ un c\^one rationnel non-singulier 
$\pi$-d\'ependant de $N^+_{\QQ}$
et soit $\sigma$ l'unique circuit inclus dans $\tau$.
{\em Dans ce qui suit, nous allons montrer que 
$X_{\tau}$
munie de l'action du sous-groupe \`a un param\`etre $\lambda _{\nu}$
est un cobordisme birationnel, que $V(\sigma)$ est l'unique
composante connexe de $\lambda _{\nu}$-points fixes de
$X_{\tau}$ et relier $\partial _-(\tau)$ et $\partial _+(\tau)$ 
\`a $(X_{\tau})_-$ et $(X_{\tau})_+$.}  

Pour cela, \'ecrivons $\tau = \langle \rho_1,\ldots, \rho_k,\ldots, 
\rho_m \rangle$
o\`u les $\rho _i$ sont les g\'en\'erateurs dans $N^+$ des faces de dimension
$1$ de $\tau$ num\'erot\'es de sorte que 
$\sigma =  \langle \rho_1,\ldots, \rho_k \rangle$. 
Comme $\tau$ est non-singulier, on peut compl\'eter la famille
$(\rho_1,\ldots, \rho_k,\ldots, \rho_m)$ en une base
$(\rho_1,\ldots,\rho_m,\ldots, \rho_{n+1})$ de $N^+$.
La vari\'et\'e torique $X_{\tau}$
est isomorphe \`a $\KK ^m \times (\KK ^*)^{n+1-m}$
et contient l'ouvert affine invariant 
$U_{\sigma} = \KK ^k \times (\KK ^*)^{n+1-k}$.
Soit $(u_1,\ldots,u_{n+1})$ la base duale 
de $(\rho_1,\ldots,\rho_{n+1})$.
Comme $\nu \in {\rm Vect} (\sigma)$, on a 
$\langle u_i,\nu\rangle =0$ pour $i \geq k+1$
et quitte \`a renum\'eroter les $\rho_j$,
on peut supposer que $\langle u_i,\nu\rangle >0$ pour $1\leq i \leq l$
et $\langle u_i,\nu\rangle <0$ pour $l+1\leq i \leq k$
o\`u $l$ est un certain entier v\'erifiant $1 \leq l \leq k-1$
(un tel $l$ existe car $\sigma$ est $\pi$-strictement convexe). 

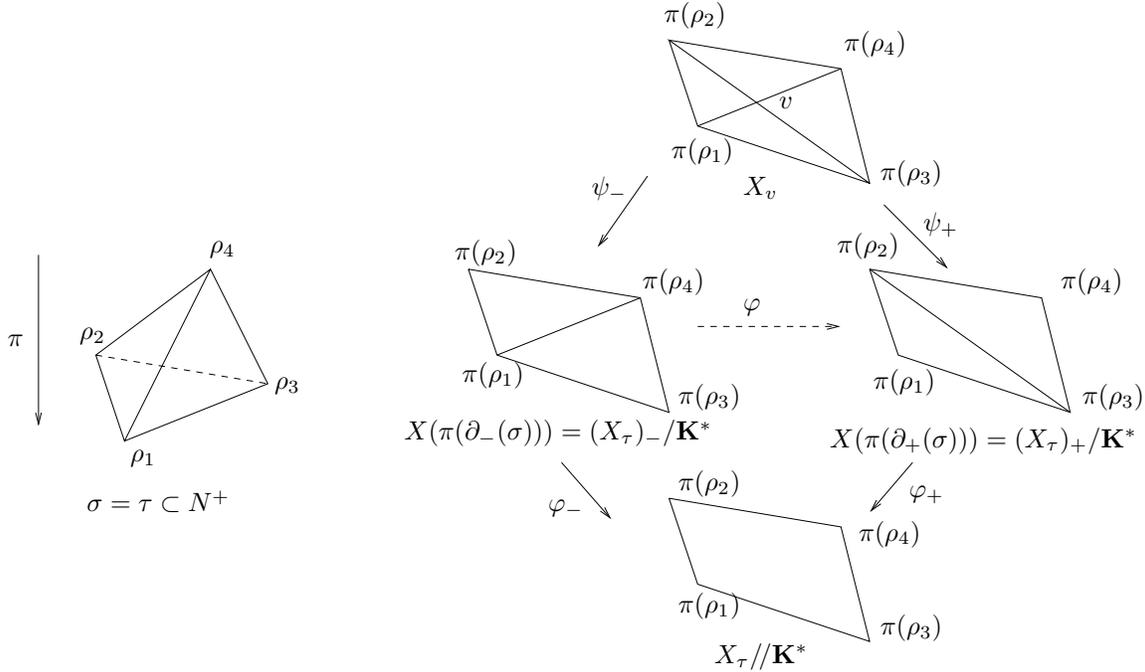
\begin{figure}[hbtp]
  \begin{center}
    \leavevmode
    \input{fig3.pstex_t}
    \caption{Cobordisme torique}
    \label{cob}
  \end{center}
\end{figure}

Le groupe $\KK ^*$ agit comme sous-groupe 
\`a un param\`etre $\lambda _{\nu}$
sur $x \in X_{\tau}$
par $$ t\cdot x = \lambda_{\nu} (t) (x_1,\ldots,x_{n+1})
= (t^{\langle u_1,\nu\rangle}x_1,\ldots, t^{\langle u_k,\nu\rangle}x_k,x_{k+1},
\dots, x_{n+1}).$$  
Pour $j$ v\'erifiant $1 \leq j \leq k$, notons
$\gamma _j = \langle \rho_1,\ldots, \check {\rho_j}, \ldots,
\rho _m \rangle $ o\`u la notation 
$\check{\rho _j}$ signifie que $\rho _j$ ne figure
pas parmi les ar\^etes de $\gamma _i$. 
C'est une face maximale de 
$\tau$ et les formules pr\'ec\'edentes montrent que :
$$ (X_{\tau})^{\KK ^*} = \{ (0,\ldots,0,x_{k+1},
\dots, x_{n+1}) \} = V(\sigma) \, ; \, 
(X_{\tau})_+ = \bigcup _{1\leq j\leq l} U_{\gamma _j} \, ; \,
(X_{\tau})_- = \bigcup _{l+1\leq j\leq k} U_{\gamma _j}.$$  
Par ailleurs, c'est un exercice facile de montrer que
$$ \partial _+ (\tau) = \bigcup _{1\leq j\leq l} \gamma _j \, \mbox{ et } 
\,
\partial _- (\tau)= \bigcup _{l+1\leq j\leq k}\gamma _j .$$   
Ainsi, la vari\'et\'e torique $X_{\pi ( \partial _+ (\tau))}$
(resp. $X_{\pi ( \partial _- (\tau))}$)
d'\'eventail
$\pi ( \partial _+ (\tau))$ (resp. $\pi ( \partial _- (\tau))$)
est la vari\'et\'e torique
$(X_{\tau})_+/\KK^*$ (resp. $(X_{\tau})_-/\KK^*$))
et $X_{\tau}$ est un cobordisme torique birationnel 
pour $ \varphi :
X_{\pi ( \partial _- (\tau))} \dra X_{\pi ( \partial _+ (\tau))}$.
Tout ceci est illustr\'e par la figure \ref{cob}. 

Etudions maintenant l'application birationnelle
$\varphi :
X_{\pi ( \partial _- (\tau))} \dra X_{\pi ( \partial _+ (\tau))}$.
Pour chaque $\rho_j$, soit 
$v_j$ le g\'en\'erateur dans $N$ de l'ar\^ete engendr\'ee
par $\pi (\rho_j)$ et $w_j$ le rationnel tel que 
$\rho_j$ soit proportionnel \`a $(v_j,w_j)$ : on \'ecrit
$\rho_j = c_j(v_j,w_j)$ o\`u $c_j$ est un entier strictement
positif.
Projetant l'\'egalit\'e 
$$ \nu = (0,1) = \sum_{j=1}^k \langle u_j, \nu \rangle \rho_j$$
sur le facteur $N_{\QQ}$ de $N_{\QQ}^+$, il vient :
$$ \sum_{j=1}^k \langle u_j, \nu \rangle c_j v_j =0.$$
Posons $$v = \sum_{j=1}^l \langle u_j, \nu \rangle c_j v_j
= - \sum_{j=l+1}^{k} \langle u_j, \nu \rangle c_j v_j \in N.$$
Ce vecteur $v$ est dans l'int\'erieur relatif 
des c\^ones $$\langle v_1,\ldots,v_l\rangle =
\bigcap _{j=l+1}^k \pi (\gamma _j) \mbox{ et }
\langle v_{l+1},\ldots,v_k\rangle =
\bigcap _{j=1}^l \pi (\gamma _j).$$ 
Par cons\'equent, les \'eventails obtenus par subdivision 
\'etoil\'ee de $\pi ( \partial _+ (\tau))$
et de $\pi ( \partial _- (\tau))$
par rapport \`a $v$ sont \'egaux, si bien qu'il existe une 
vari\'et\'e torique $X_v$ et un diagramme commutatif~:

\centerline{
\xymatrix{ & X_v \ar[ld]_{\psi_-}   \ar[rd]^{\psi _+}&\\
(X_{\tau})_-/\KK^* \ar[rd]_{\varphi_-} \ar@{-->}[rr]^{\varphi} & & 
(X_{\tau})_+/\KK^* \ar[ld]^{\varphi_+}\\
& X_{\tau}/\!/\KK^* & 
}
}
Si de plus le c\^one initial $\tau$ est $\pi$-non 
singulier, alors $(X_{\tau})_-/\KK^*$ et $(X_{\tau})_+/\KK^*$
sont lisses et
$$v = \sum_{j=1}^l v_j
=  \sum_{j=l+1}^{k} v_j \in N,$$
de sorte que $X_v$ est lisse, 
$\psi_+ : X_v \to (X_{\tau})_+/\KK^*$ et 
$\psi_- : X_v \to (X_{\tau})_-/\KK^*$ sont des \'eclatements 
le long de sous-vari\'et\'es lisses invariantes (voir la figure
\ref{cob}).

\subsection{Cobordisme torique}

Soit $\varphi : X_{\Sigma '}  \to X_{\Sigma}$ 
une application r\'eguli\`ere
birationnelle \'equivariante entre deux vari\'et\'es toriques
projectives lisses de dimension $n$
et d'\'eventails respectifs $\Sigma '$ et $\Sigma$
dans $N_{\QQ}$. 
La construction expliqu\'ee au paragraphe \ref{existe cobord}
s'adapte sans difficult\'e
au cas torique (historiquement, rappelons que
c'est la construction torique qui a inspir\'e 
le cas g\'en\'eral), si bien qu'il existe
un cobordisme birationnel $B$ torique (au sens o\`u l'action de
$\KK ^*$ sur $B$ est celle d'un sous-groupe \`a un param\`etre
du tore de $B$)
quasi-projectif
pour $\varphi$.

D'apr\`es ce qui pr\'ec\`ede, l'existence d'un tel $B$ 
en termes d'\'eventails se traduit de la fa\c con
suivante. Il existe un \'eventail 
$\Sigma$ non-singulier dans $N_{\QQ}^+$
tel que $\pi (\partial _- (\Sigma)) = \Sigma '$,
$\pi (\partial _+ (\Sigma)) = \Sigma$, $B$ est la vari\'et\'e 
torique d'\'eventail $\Sigma$ et $\varphi$
est le cobordisme birationnel $\varphi : B_-/\KK^* = X_{\Sigma'}
\to B_+/\KK^* = X_{\Sigma}$.
Comme $B$
est quasi-projectif, il y a un ordre sur
les circuits de $\Sigma$ (qui correspondent
aux composantes connexes des points fixes
pour l'action du sous-groupe \`a un param\`etre $\nu = (0,1)$).
Choisissons un circuit $\sigma_1$ minimal pour cet ordre. 
Soient 
${\rm Star}(\sigma_1) = \{ \tau \in \Sigma \, |\,
\sigma_1 \subset \tau \}$
et $\overline{{\rm Star}}(\sigma_1) = \{ \tau' \in \Sigma \, |\,
\tau' \subset \tau  \mbox{ pour un } \tau \mbox{ dans } 
{\rm Star}(\sigma_1) \}.$  
Alors comme $\sigma_1$
est minimal, $\partial _-(\overline{{\rm Star}}(\sigma_1)) 
\subset \partial _-(\Sigma)$. Posons
alors
$\Sigma _1 = (\Sigma \setminus \partial _-(\overline{{\rm Star}}(\sigma_1)))
\cup \partial _+(\overline{{\rm Star}}(\sigma_1))$.
D'apr\`es le paragraphe pr\'ec\'edent, si $B_1$
est la vari\'et\'e d'\'eventail $\Sigma_1$ et si 
$X_1$ est la vari\'et\'e torique d'\'eventail
$\Sigma _1 :=\pi(\partial _-(\Sigma _1))$, alors $\varphi$ se
d\'ecompose en
$$\varphi : B_-/\KK^* = X_{\Sigma '} \stackrel{\psi_-}{\leftarrow}
X_v \stackrel{\psi_+ }{\rightarrow} X_{\Sigma _1}  = 
(B_1)_-/\KK^*  \stackrel{\varphi_1}{\dra} X_{\Sigma}= (B_1)_+/\KK^*
= B_+/\KK^* ,$$
et on peut recommencer avec le cobordisme $B_1$.
Si de plus $\Sigma$ est $\pi$-non singulier, alors 
$X_{\Sigma _1}$ et $X_v$ sont lisses, 
$\psi_-$ et $\psi_+ $ sont des \'eclatements 
le long de sous-vari\'et\'es lisses invariantes
et $\Sigma_1$ est encore $\pi$-non singulier.

\medskip

{\bf Faisons le point :} le th\'eor\`eme \ref{factorisation torique}
est d\'emontr\'e pour une application r\'eguli\`ere 
$\varphi : X_{\Sigma '}  \to X_{\Sigma}$
birationnelle \'equivariante entre deux vari\'et\'es toriques
projectives lisses de dimension $n$
si on peut contruire un cobordisme 
birationnel $\tilde B$ quasi-projectif 
d'\'eventail $\tilde{\Sigma}$ {\bf $\pi$-non singulier}
dans $N_{\QQ}^+$ tel que 
$\pi (\partial _- (\tilde{\Sigma})) = \Sigma '$ et
$\pi (\partial _+ (\tilde{\Sigma})) = \Sigma$.

Un tel \'eventail $\tilde{\Sigma}$ va \^etre obtenu
en $\pi$-d\'esingularisant l'\'eventail $\Sigma$, son existence est 
garantie par le th\'eor\`eme de $\pi$-d\'esingularisation de Morelli,
qui fait l'objet du paragraphe suivant.

\subsection{Le th\'eor\`eme de $\pi$-d\'esingularisation}

Les notations sont celles du paragraphe pr\'ec\'edent.
Le th\'eor\`eme suivant
est d\^u \`a Morelli \cite{Mor96} avec une d\'emonstration 
incompl\`ete~; la d\'emonstration a \'et\'e compl\'et\'ee 
par Abra\-mo\-vich, Matsuki
et Rashid \cite{AMR99} \cite{Mat00}. 

\begin{theo}\label{pi desingularisation}
Soit $\Sigma$ un \'eventail de $N^+_{\QQ}$. Alors il existe un \'eventail
$\pi$-non-singulier
de $N^+_{\QQ}$ obtenu par une suite 
finie de subdivisions
\'etoil\'ees de $\Sigma$. De plus, 
la proc\'edure de $\pi$-d\'esingula\-risa\-tion n'affecte pas 
les c\^ones $\pi$-non-singuliers de $\Sigma$.
\end{theo}
 
\subsection{D\'emonstration du th\'eor\`eme de factorisation faible torique}
On termine la d\'emonstration du th\'eor\`eme de factorisation, qui
repose sur le lemme de Chow torique (voir par exemple \cite{Oda88})
et la lev\'ee des ind\'eterminations
torique due \`a De Concini-Procesi \cite{DCP85}.

\begin{theo}\label{chow torique} 
Soit $X$ une vari\'et\'e torique compl\`ete.
Alors il existe une vari\'et\'e torique projective
lisse $\tilde X$ et une application r\'eguli\`ere birationnelle
\'equivariante $\varphi : \tilde X \to X$. 
\end{theo}

\begin{theo}\label{concini}
Soit $\varphi : X \dra X'$ une application birationnelle
\'equivariante entre deux vari\'et\'es toriques compl\`etes lisses.
Alors il existe une 
suite d'\'eclatements le long de sous-vari\'et\'es
inva\-rian\-tes de codimension deux $\psi :\bar X \to X$ 
telle que $\varphi \circ \psi : \bar X \to X'$ soit 
une application r\'eguli\`ere birationnelle
\'equivariante.
\end{theo}

De ces deux r\'esultats, d\'eduisons une version pr\'ecise du
th\'eor\`eme de Moishezon torique :

\begin{theo}\label{moishezon torique}
Soit $X$ une vari\'et\'e torique lisse et compl\`ete.
Alors il existe une vari\'et\'e torique projective
lisse $\tilde X$ obtenue \`a partir de $X$
par suite d'\'eclatements le long de sous-vari\'et\'es
invariantes de codimension deux $\psi : \tilde X \to X$.
\end{theo}

\noindent {\em D\'emonstration.} 
Soit $\varphi : \tilde X \to X$ donn\'ee par le 
th\'eor\`eme \ref{chow torique}, o\`u $\tilde X$ est projective lisse. 
Appliquons le th\'eor\`eme \ref{concini} \`a $\varphi ^{-1}$ : il 
y a une suite d'\'eclatements le long de sous-vari\'et\'es
invariantes de codimension deux $\psi : \bar X \to X$
de sorte que $h=\varphi ^{-1} \circ \psi : \bar X \to \tilde X$
soit une application r\'eguli\`ere birationnelle.
Comme $\psi$ est un morphisme projectif et puisque
$\varphi \circ h = \psi$, $h$ est aussi
un morphisme projectif, et comme $\tilde X$ est projective,
on en d\'eduit que $\bar X$ est projective. \finpreuve

\medskip

{\bf Application :} soit $\varphi : X_1 \dra X_2$
une application birationnelle entre deux vari\'et\'es 
toriques compl\`etes et lisses. 
Par le th\'eor\`eme \ref{moishezon torique}, 
il y a des vari\'et\'es toriques projectives lisses $\tilde{X_1}$
et $\tilde{X_2}$ et deux suites d'\'eclatements le long de sous-vari\'et\'es
invariantes de codimension deux $\psi_1 : \tilde{X_1}\to X_1$
et $\psi_2 : \tilde{X_2}\to X_2$.
Soit $\tilde{\varphi} = \psi_2^{-1}\circ \varphi \circ \psi_1 :
\tilde{X_1} \dra \tilde{X_2}$.
Par le th\'eor\`eme \ref{concini}, il existe
une 
suite d'\'eclatements le long de sous-vari\'et\'es
invariantes de codimension deux $\psi :\bar {X_1} \to \tilde{X_1}$ 
telle que $\tilde{\varphi} \circ \psi : \bar{X_1} \to \tilde{X_2}$ soit 
une application r\'eguli\`ere birationnelle
\'equivariante. Comme $\bar{X_1}$ et $\tilde{X_2}$ sont projectives lisses,
il y a d'apr\`es les paragraphes pr\'ec\'edents une factorisation
faible pour $\tilde{\varphi} \circ \psi$ et par suite
une factorisation faible pour $\varphi$.
Ceci ach\`eve la d\'emonstration du th\'eor\`eme 
\ref{factorisation torique}. 

\subsection{Extension au cas toro\"{\i}dal}
Il y a une classe de vari\'et\'es alg\'ebriques poss\'edant
essentiellement les m\^emes propri\'et\'es que les
vari\'et\'es toriques. 

\medskip

\noindent {\bf D\'efinition.} Soit $V$ une vari\'et\'e
alg\'ebrique et $U$ un ouvert de $V$.
Le plongement $U \subset V$ est un {\em plongement
toro\"{\i}dal} si pour tout $x \in V$, il y a un voisinage
$V_x$ de $x$ dans $V$, une vari\'et\'e torique $X_x$
et une application \'etale $\eta_x : V_x \to X_x$
telle que $ \eta_x ^{-1} (T_x) = V_x \cap U$. 
 
\medskip

Comme dans le cas torique, il y a une notion
naturelle de morphismes toro\"{\i}daux  
et d'ap\-pli\-ca\-tions bi\-ration\-nel\-les toro\"{\i}dales 
entre plongements
toro\"{\i}daux. De plus, il y a une stratification 
naturelle de $V$ et les adh\'erences des strates
correspondent formellement aux adh\'erences d'orbites dans
les vari\'et\'es toriques.
Enfin, un plongement toro\"{\i}dal est d\'ecrit
par un complexe poly\'edral conique rationnel qui
joue exactement le r\^ole de l'\'eventail d'une vari\'et\'e 
torique (\`a ceci pr\`es que ce complexe est d\'efini abstraitement et
qu'il ne se plonge en g\'en\'eral pas lin\'eairement dans un espace vectoriel
$N_{\QQ}$). Les r\'esultats de la g\'eom\'etrie birationnelle
des vari\'et\'es toriques (existence de d\'esingularisation, 
lev\'ee des ind\'eterminations par \'eclatements de centres lisses, etc.)
s'\'etendent sans trop de difficult\'es au cas des plongements
toro\"{\i}daux (voir \cite{KKMS73} ou 
\cite{AbO97} p. 69-73).
Il en est de m\^eme pour le th\'eor\`eme de factorisation comme
l'ont remarqu\'e Abramovich, Matsuki et Rashid \cite{AMR99}.

\begin{theo}\label{factorisation toroidale}
Soit $\varphi: (U_1 \subset X_1) \dashrightarrow (U_2 \subset X_2)$ 
une application birationnelle toro\"{\i}dale entre deux 
plongements toro\"{\i}daux lisses sur $\KK$, $X_1$ et $X_2$ compl\`etes. 
Alors, $\varphi$ se factorise en 
une suite d'\'eclatements et de contractions de centres 
lisses \'egaux \`a des adh\'erences de strates.
\end{theo} 
 
\section{Structures toro\"{\i}dales et 
cobordisme birationnel}\label{struc toroid}
Le but de cette partie est 
de d\'ecomposer un cobordisme birationnel donn\'e
en une suite de cobordismes birationnels (en g\'en\'eral
singuliers) 
toro\"{\i}daux. L'id\'ee est la suivante : si on \'eclate un faisceau
d'id\'eaux $I$ sur une vari\'et\'e alg\'ebrique, on obtient une nouvelle
vari\'et\'e alg\'ebrique avec un diviseur bien d\'etermin\'e
puisque $I$ devient principal apr\`es \'eclatement. 
Le compl\'ementaire de son support est alors un candidat pour 
\^etre un plongement toro\"{\i}dal. Lorsque de plus 
la vari\'et\'e alg\'ebrique est munie d'une action de $\KK^*$,
on peut essayer de construire $I$ de sorte que 
le plongement toro\"{\i}dal obtenu soit compatible avec l'action de $\KK^*$.

\subsection{Id\'eal $\alpha$-toro\"{\i}dal}
Soit $V$ une vari\'et\'e alg\'ebrique munie d'une action 
alg\'ebrique effective de $\KK ^*$. Pour $v \in V$, on note
$\Stab (v)$ le stabilisateur de $v$ 
($\Stab (v)$ est isomorphe soit \`a $\KK^*$,
soit \`a un groupe cyclique $\ZZ / n\ZZ$). 

\medskip

\noindent {\bf Notation.} Si $\alpha$ est un entier,
on note $J_{\alpha,v}$ l'id\'eal de ${\mathcal O}_{V,v}$
engendr\'e par les fonctions $f \in {\mathcal O}_{V,v}$
$\Stab (v)$-semi-invariantes de $\Stab (v)$-poids $\alpha$,
{\em i.e.} pour tout $t \in \Stab (v)$, 
$t^*(f) = t^{\alpha}f$. 
 
\medskip

Il n'est pas difficile de voir que si
$z_1,\ldots,z_n$ sont des g\'en\'erateurs $\Stab (v)$-semi-invariants
de l'id\'eal maximal ${\mathcal M}_{V,v}$ dont l'existence
est assur\'ee par le th\'eor\`eme de Sumihiro \cite{Sum74} \cite{Sum75},
alors $J_{\alpha,v}$ est engendr\'e par les mon\^omes
$z_1^{m_1}\cdots z_n ^{m_n}$ 
de $\Stab (v)$-poids $\alpha$ ({\em i.e.} v\'erifiant
$a_1 m_1+\cdots+a_n m_n =\alpha$ si $a_i$ d\'esigne
le $\Stab (v)$-poids de $z_i$). Ceci implique en particulier
que si $\eta : V \to X$ est 
un morphisme $\KK ^*$-\'equivariant fortement \'etale 
entre deux vari\'et\'es alg\'ebriques munies
d'une action de $\KK ^*$, alors
l'image inverse de $J_{\alpha,\eta(v)}$ engendre
$J_{\alpha,v}$ pour tous $v \in V$ et $\alpha \in \ZZ$.

La collection des $J_{\alpha,v}$ lorsque $v$ d\'ecrit
$V$ ne d\'efinit pas un faisceau d'id\'eaux coh\'erent en g\'en\'eral.
Cependant, on a la proposition suivante :

\begin{prop}\label{ideal toroidal} 
Soient $\alpha \in \ZZ$ et $B$ un cobordisme birationnel
lisse  
quasi-projectif de la forme $B_{a_i}$ 
(les notations sont celles des parties \ref{GIT} et \ref{slt}).
Alors, il existe un unique faisceau d'id\'eaux coh\'erent
$\KK ^*$-\'equivariant $I_{\alpha}$ non nul tel que pour tout
$v \in (B_+\cap B_-)\cup B^{\KK ^*}$, on a 
$(I_{\alpha})_v = J_{\alpha,v}$.
\end{prop}

L'id\'eal $I_{\alpha}$ s'appelle le faisceau {\em $\alpha$-toro\"{\i}dal}
associ\'e \`a l'action de $\KK ^*$ sur $B$.
Rappelons que $(B_+\cap B_-)$ est un ouvert non vide de $B$
et que $(B_+\cap B_-)\cup B^{\KK ^*}$ est l'union des $\KK ^*$-orbites
ferm\'ees de $B$.

\medskip

\noindent {\em D\'emonstration de la proposition 
\ref{ideal toroidal}.} 
Commen\c cons par l'unicit\'e de $I_{\alpha}$. 
Elle est claire sur $(B_+\cap B_-)\cup B^{\KK ^*}$. 
Soit donc $v$ dans le compl\'ementaire de $(B_+\cap B_-)\cup B^{\KK ^*}$.
L'adh\'erence de l'orbite de $v$ contient un point fixe.
Comme $I_{\alpha}$ est 
uniquement d\'etermin\'e en ce point fixe, il l'est aussi au voisinage
puisque $I_{\alpha}$ est coh\'erent. Par $\KK ^*$-\'equivariance,
$I_{\alpha}$ est alors
uniquement d\'etermin\'e en $v$.

Pour l'existence, la proposition \ref{cartes toriques}
permet de recouvrir $B$ par des cartes toriques 
$\eta_x : V_x \to X_x$ fortement
\'etales. Comme $((B_+\cap B_-)\cup B^{\KK ^*}) \cap V_x$
est l'image inverse par $\eta_x$ 
de $((X_x)_+\cap (X_x)_-)\cup (X_x)^{\KK ^*}$,
il suffit de construire un faisceau 
$\alpha$-toro\"{\i}dal  
sur $X_x$. 
Par image inverse, ceci d\'efinit un faisceau 
$\alpha$-toro\"{\i}dal 
sur $V_x$ et l'unicit\'e d\'ej\`a prouv\'ee entra\^{\i}ne
que ces faisceaux ainsi construits se recollent en un
faisceau sur $B$. Le lemme suivant \cite{AKMW99} ach\`eve donc
la preuve.
\finpreuve   
 
\begin{lemm}\label{exist torique}
Soit $X$ une vari\'et\'e torique affine lisse
d\'efinie par un c\^one r\'egulier $\sigma$ d'un r\'eseau $N$
et soit $a\in N$ un \'el\'ement primitif correspondant 
\`a une action effective de $\KK^*$ 
d'un sous-groupe \`a un param\`etre
de $T$. On suppose de plus que $a$ n'appartient
pas \`a $\sigma \cup -\sigma$ (ceci correspond au fait
que $B_+ \cap B_-$ est non vide). 
Alors, pour tout $\alpha \in \ZZ$, le faisceau $I_{\alpha}$ existe,
est non nul et engendr\'e par les mon\^omes
$z^m$ o\`u $m\in \check {\sigma}$ avec $\langle m,a \rangle=\alpha$.
\end{lemm} 

Remarquons qu'il n'est pas utile d'apr\`es ce lemme
de calculer explicitement les stabilisateurs $\Stab (v)$. 
Donnons un exemple concret, que nous utiliserons aussi 
dans la suite : 

\medskip

\noindent {\bf Exemple ``$(2,1,-1)$''.} Soit $X = \KK ^3$
munie de l'action du sous-groupe \`a un 
param\`etre $a =(2,1,-1)$ ({\em i.e.} $\lambda _a(t)(z_1,z_2,z_3)
= (t^2 z_1,t z_2,t^{-1} z_3)$), alors $I_2$ (resp. $I_1$, resp. $I_{-1}$)
est le faisceau d'id\'eaux  
engendr\'e par $z_1$ et $z_2 ^2$
(resp. $z_1 ^2 z_3$ et $z_2$, resp. $z_3$).

\subsection{Plongements toro\"{\i}daux avec action toro\"{\i}dale
de $\KK ^*$}
Soit $V$ une vari\'et\'e alg\'ebrique munie d'une action 
alg\'ebrique de $\KK ^*$. On suppose que l'action de $\KK ^*$
est localement torique au sens o\`u
on peut recouvrir $V$ par des cartes toriques 
fortement \'etales.
Soit $D$ un diviseur effectif de $V$
et $U$ l'ouvert $V \setminus \Supp(D)$. 

\medskip

\noindent {\bf D\'efinition.}
On dit que $U$ est un {\em plongement  
toro\"{\i}dal avec action toro\"{\i}dale de $\KK^*$}
si pour tout $x \in V$, il y a un voisinage  
$\KK^*$-invariant et affine $V_x$ de $x$, 
une vari\'et\'e torique affine $X_x$ avec une action
de $\KK ^*$ d'un sous-groupe \`a un param\`etre 
et un morphisme $\KK ^*$-\'equivariant \'etale 
$\eta_x : V_x \to X_x$
tel que $\eta_x^{-1} (T_x)= V_x \cap U$
(o\`u $T_x$ d\'esigne le tore dense de $X_x$).
On appelle {\em carte toro\"{\i}dale \'etale} en $x$ la donn\'ee
d'une telle $\eta_x : V_x \to X_x$. D'apr\`es le    
lemme fondamental de Luna, si une telle carte existe,
en restreignant l'ouvert $V_x$, on peut la supposer 
de plus fortement \'etale.
 
\medskip
 
\noindent {\bf Exemples.} 
Soit $X$ la vari\'et\'e torique 
$\KK ^2$ d\'efinie par un c\^one r\'egulier $\sigma = \langle v_1,v_2\rangle$ 
d'un r\'eseau $N$ et soit $a =(a_1,a_2) \in N$ 
un \'el\'ement primitif correspondant 
\`a une action effective de $\KK^*$ 
d'un sous-groupe \`a un param\`etre
de $T$. Soit $D_{v_1} = \{ (z_1,z_2) \in \KK ^2 \, | \, z_1=0 \}$ 
et $D_{v_2} = \{ (z_1,z_2) \in \KK ^2 \, | \, z_2=0 \}$.
Alors, $U = \KK ^2 \setminus (D_{v_1} \cup D_{v_2})$ 
est un plongement  
toro\"{\i}dal avec action toro\"{\i}dale de $\KK^*$ quel que soit
le choix du sous-groupe \`a un param\`etre $a$.
En revanche, $U = \KK ^2 \setminus D_{v_1}$ 
est un plongement  
toro\"{\i}dal avec action toro\"{\i}dale de $\KK^*$
si et seulement si $a = \pm v_1$. 
G\'eom\'etriquement, ceci signifie que l'action de
$\KK^*$ sur $\KK ^2$
est produit de l'action triviale sur $\{0\} \times \KK$
et d'une action d'un sous-groupe \`a un param\`etre
de la vari\'et\'e torique $\KK \times \{0\}$ 
et que le diviseur $D_{v_2}$ est \'egal \`a $\KK \times \{0\}$.

\medskip

Cet exemple se g\'en\'eralise en le lemme suivant \cite{AKMW99} :

\begin{lemm}\label{caract torique}  
Soit $X$ une vari\'et\'e torique 
d\'efinie par un \'eventail $\Sigma$ 
d'un r\'eseau $N$ et soit $a \in N$ 
un \'el\'ement primitif correspondant 
\`a une action effective de $\KK^*$ 
d'un sous-groupe \`a un param\`etre
de $T$, soient $D$ un diviseur torique effectif de $X$
et $U = X \setminus \Supp(D)$.
Les assertions suivantes sont \'equivalentes :
\begin{enumerate}
\item [(i)] 
pour tout c\^one $\sigma$ de $\Sigma$
et tout diviseur torique $E$ de l'ouvert 
torique affine $U_{\sigma} \subset X$, 
si $E$ n'est pas dans $\Supp(D)$, alors il existe une 
vari\'et\'e torique affine $X_{\sigma '}$
telle que :
$$ U_{\sigma} \simeq X_{\sigma '} \times \KK \, \mbox{ et }
\, E \simeq X_{\sigma '} \times \{0\} $$
de sorte que l'action de $\KK^*$ sur $ U_{\sigma}$
soit le produit de l'action d'un sous-groupe \`a un param\`etre
sur $X_{\sigma '}$ et de l'action triviale sur $\KK$,
\item[(ii)] $U$ est un plongement  
toro\"{\i}dal avec action toro\"{\i}dale de $\KK^*$.
\end{enumerate}
\end{lemm}

Remarquons que l'implication {\em (i) implique (ii)}
est ais\'ee, c'est la seule utilis\'ee dans la suite.

\subsection{Faisceaux toro\"{\i}daux et plongements toro\"{\i}daux}
Soit $B$ un cobordisme birationnel
lisse et quasi-projectif de dimension $n+1$ 
de la forme $B_{a_i}$
(les notations sont celles des parties \ref{GIT} et \ref{slt}).
Fixons nous une famille finie de
cartes toriques fortement \'etales $\eta_x : V_x \to X_x$
recouvrant $B$. Dans $X_x$, l'action de $\KK^*$ correspond
\`a celle d'un sous-groupe \`a un param\`etre : concr\`etement
$X_x$ est isomorphe \`a $\KK ^m \times (\KK ^*)^{n+1 -m}$ 
et l'action de $\KK^*$ est de la forme
$t \cdot (x_1,\ldots,x_{n+1}) = (t^{\alpha_{1}}x_1,\ldots,
t^{\alpha_{n+1}}x_{n+1})$ pour certains poids entiers 
$\alpha_1,\ldots,\alpha_{n+1}$. 

\medskip

\noindent {\bf Notation.} On note ${\mathcal A}$
une famille finie d'entiers telle que pour toute carte
torique fortement \'etale $\eta_x : V_x \to X_x$,
tous les poids $\alpha_1,\ldots,\alpha_{n+1}$ 
de l'action de $\KK ^*$ appartiennent
\`a ${\mathcal A}$.
Une telle famille d'entiers sera dite {\em admissible}.

\medskip

\noindent {\bf D\'efinition.}
Si ${\mathcal A}$ est une famille admissible d'entiers,
on note $I_{\mathcal A}$ le faisceau d'id\'eaux sur $B$ d\'efini
par $I_{\mathcal A} = \prod_{\alpha \in {\mathcal A}} I_{\alpha}$.
Un tel faisceau d'id\'eaux sur $B$ est appel\'e {\em faisceau
toro\"{\i}dal}.
  
\medskip

\noindent {\bf Notation.}
Si $I_{\mathcal A}$ est un faisceau
toro\"{\i}dal, on note $\pi_{\mathcal A} : B_{I_{\mathcal A}}^{tor} \to B$
la normalisation de l'\'eclatement de $B$ 
de centre le faisceau d'id\'eaux $I_{\mathcal A}$. 
Comme $I_{\mathcal A}$ est $\KK ^*$-\'equivariant, 
la vari\'et\'e $B_{I_{\mathcal A}}^{tor}$ est naturellement
munie d'une action alg\'ebrique de $\KK ^*$ faisant
de $\pi_{\mathcal A} $ un morphisme \'equivariant.
On note $D_{I_{\mathcal A}}^{tor}$ le diviseur 
d\'efini par le faisceau inversible 
$\pi_{\mathcal A}^{-1}(I_{\mathcal A})
{\mathcal O}_{B_{I_{\mathcal A}}^{tor}}$ et 
$U_{I_{\mathcal A}}^{tor}$
le compl\'ementaire de son support dans $B_{I_{\mathcal A}}^{tor}$.  

La proposition suivante est le coeur de cette partie :

\begin{prop}
La vari\'et\'e alg\'ebrique $B_{I_{\mathcal A}}^{tor}$ est un 
cobordisme birationnel quasi-projectif. 
L'ouvert $U_{I_{\mathcal A}}^{tor}$ est un 
plongement toro\"{\i}dal avec action toro\"{\i}dale de $\KK^*$.
Enfin, $(B_{I_{\mathcal A}}^{tor})_+ = \pi_{\mathcal A}^{-1}(B_+)$
et $(B_{I_{\mathcal A}}^{tor})_- = \pi_{\mathcal A}^{-1}(B_-)$.
\end{prop}

\noindent {\em D\'emonstration.} 
Montrons que l'ouvert $U_{I_{\mathcal A}}^{tor}$ est un 
plongement toro\"{\i}dal avec action toro\"{\i}dale de $\KK^*$.
A nouveau, 
il suffit de traiter le cas torique, autrement dit de montrer
ce r\'esultat dans chacune des cartes toriques fortement \'etales 
consid\'er\'ees au d\'ebut de la construction. Le lemme suivant
donne donc le r\'esultat. \finpreuve 

\begin{lemm}
Soit $X$ la vari\'et\'e torique affine lisse 
$\KK ^m \times (\KK ^*)^{n+1 -m}$ 
avec une action de $\KK^*$ de la forme
$t \cdot (x_1,\ldots,x_{n+1}) = (t^{\alpha_{1}}x_1,\ldots,
t^{\alpha_{n+1}}x_{n+1})$ pour certains poids entiers 
$\alpha_1,\ldots,\alpha_{n+1}$. 
Soit ${\mathcal A}$ une famille d'entiers 
contenant $\alpha_1,\ldots,\alpha_{n+1}$ et 
$I_{\mathcal A}$ le faisceau toro\"{\i}dal associ\'e.
Soient $\pi_{\mathcal A} : X_{I_{\mathcal A}}^{tor}\to X$
la normalis\'ee de l'\'eclatement de centre $I_{\mathcal A}$,
$D_{I_{\mathcal A}}^{tor}$ le diviseur
d\'efini par le faisceau inversible 
$\pi_{\mathcal A}^{-1}(I_{\mathcal A})
{\mathcal O}_{B_{I_{\mathcal A}}^{tor}}$
et 
$U_{I_{\mathcal A}}^{tor}$
le compl\'ementaire de son support dans $B_{I_{\mathcal A}}^{tor}$.
Alors $U_{I_{\mathcal A}}^{tor}$ est un 
plongement toro\"{\i}dal avec action toro\"{\i}dale de $\KK^*$.
\end{lemm}

\noindent {\em D\'emonstration.} 
Les diviseurs toriques de $X_{I_{\mathcal A}}^{tor}$ non 
inclus dans $\Supp(D_{I_{\mathcal A}}^{tor})$
(les seuls \`a consid\'erer d'apr\`es le lemme
\ref{caract torique}) sont des transform\'ees strictes
de diviseur $D_i := \{ x \in X \, |\, x_i=0\}$.
Soit donc $i$ dans $\{1,\ldots, n+1\}$. Evidemment,
le mon\^ome $x_i$ appartient au 
faisceau $\alpha_i$-toro\"{\i}dal $I _{\alpha_i}$.
Si ce dernier est principal engendr\'e par $x_i$, alors 
$D_i$ est dans $\Supp (D_{I_{\mathcal A}}^{tor})$. 
Sinon, soit $X_i$ la normalis\'ee de 
l'\'eclatement du faisceau $\alpha_i$-toro\"{\i}dal $I _{\alpha_i}$.
Il est ais\'e de v\'erifier que $D_i$
satisfait l'assertion (i) 
du lemme \ref{caract torique}. La suite de la d\'emonstration   
consiste \`a remarquer que $X_{I_{\mathcal A}}^{tor}$
est obtenu en normalisant l'\'eclatement de centre $I_{\alpha_1}$,
puis en normalisant l'\'eclatement de centre l'image inverse
de $I_{\alpha_2}$, etc.\finpreuve

\medskip

Reprenons l'exemple ``$(2,1,-1)$'' : $X = \KK ^3$
munie de l'action du sous-groupe \`a un param\`etre
$(2,1,-1)$.
On prend ici ${\mathcal A} = \{ -1,1,2\}$.
Si $X$ est d\'efinie par un c\^one r\'egulier $\sigma = 
\langle v_1,v_2,v_3 \rangle$, alors 
$X_{I_{\mathcal A}}^{tor}$ est la vari\'et\'e torique d'\'eventail
$\Sigma$ obtenu en subdivisant $\sigma$ en les trois c\^ones :
$\sigma_1 = \langle v_1,2v_1+v_2,v_3\rangle$,
$\sigma_2 = \langle v_1+v_2,v_3,2v_1+v_2,v_2+v_3 \rangle$
et $\sigma_3 = \langle v_1+v_2,v_2+v_3,v_2 \rangle$.
Le diviseur $D_{I_{\mathcal A}}^{tor}$ est le diviseur
$D_{v_3} \cup D_{v_1+v_2} \cup D_{2v_1+v_2} \cup D_{v_2+v_3}$.
Seul les c\^ones $\sigma_1$ et $\sigma_3$ sont
\`a consid\'erer pour v\'erifier
l'assertion (i) du lemme \ref{caract torique}, avec
$E = D_{v_1}$ pour $\sigma_1$ et $E= D_{v_2}$
pour $\sigma_3$. L'assertion (i) est satisfaite car 
$a=(2,1,-1)$ appartient au r\'eseau engendr\'e par 
$v_1+v_2$ et $v_2+v_3$
et au r\'eseau engendr\'e par $2v_1+v_2$ et $v_3$.

\begin{coro}\label{coro toro}
Soit $B_{a_i}$ un cobordisme birationnel lisse et quasi-projectif obtenu
dans le paragraphe \ref{GIT}. 
Choisissons un faisceau toro\"{\i}dal $I_{{\mathcal A}_i}$
et notons $\pi _{{\mathcal A}_i} 
: B_{a_i}^{tor} \to B_{a_i}$ la normalis\'ee de
l'\'ecla\-te\-ment de centre $I_{{\mathcal A}_i}$.
Alors, il y a un diagramme commutatif :

\centerline{ 
\xymatrix{ & 
(B_{a_i}^{tor})_-/\KK^* \ar[dl]_{(\pi_{{\mathcal A}_i})_-} 
\ar@{-->}[r]^{\varphi_i^{tor}} &
(B_{a_i}^{tor})_+/\KK^* \ar[dr]^{(\pi_{{\mathcal A}_i})_+} \\
V_i = (B_{a_i})_-/\KK^*   \ar@{-->}[rrr]^{\varphi_i} & &  & 
(B_{a_i})_+/\KK^* = V_{i+1} 
}
} 
\end{coro}

\medskip

{\bf Faisons le point :} 
la vari\'et\'e $B_{a_i}^{tor}$ contient l'ouvert 
$U_{I_{{\mathcal A}_i}}^{tor}$ comme plongement toro\"{\i}dal avec action de
$\KK^*$ si bien que $\varphi_i^{tor}$ est une application birationnelle
toro\"{\i}dale au sens o\`u pour tout $x \in B_{a_i}^{tor}$
il y a une carte en $x$ toro\"{\i}dale fortement \'etale $\eta_x : V_x^{tor} 
\to X_x$ (donc v\'erifiant $\eta_x ^{-1}(T_x) =
V_x^{tor}\cap U_{I_{{\mathcal A}_i}}^{tor}$)
induisant un diagramme commutatif : 

\centerline{ 
\xymatrix{  
(V_x ^{tor})_-/\KK^* \ar[d]_{(\eta_x)_-/\KK^*} 
\ar@{-->}[r]^{\varphi_i^{tor}} &
(V_x^{tor})_+/\KK^* \ar[d]^{(\eta_x)_+/\KK^*} \\
(X_x)_-/\KK^*  \ar@{-->}[r] & (X_x)_+/\KK^* 
}
}  
\noindent o\`u $(X_x)_-/\KK^*  \dra (X_x)_+/\KK^*$ est une 
application birationnelle
entre vari\'et\'es toriques. Dans la partie suivante, on explique
comment on peut obtenir un r\'esultat analogue en restant
dans le cadre des vari\'et\'es non singuli\`eres.

\section{D\'esingularisation canonique et d\'emonstration du
th\'eor\`eme de factorisation}\label{desing}

\subsection{D\'esingularisation canonique}
Sur un corps alg\'ebriquement clos de caract\'eristique nulle,
on sait depuis Hironaka \cite{Hir64} que toute vari\'et\'e alg\'ebrique
peut \^etre d\'esingularis\'ee par une suite d'\'eclatements
le long de centres lisses. A partir de la dimension trois,
il n'y a pas de d\'esingularisation minimale
et il n'y a pas de choix naturel
d'un mod\`ele non singulier pour une vari\'et\'e alg\'ebrique
donn\'ee. Cependant, avec les travaux de Bierstone-Milman
\cite{BiM97},
Encinas-Villamayor \cite{EnV97} et Villamayor \cite{Vil89}, 
on peut parler   
de d\'esingularisation ``canonique''.

\medskip

\noindent {\bf D\'efinition.} Une 
{\em r\'esolution canonique des singularit\'es} 
est un algorithme qui \`a toute vari\'et\'e alg\'ebrique $X$
associe une suite uniquement d\'etermin\'ee d'\'eclatements
le long de centres lisses $r: X^{res} \to X$ satisfaisant :
pour tout morphisme lisse $Y\to X$, la vari\'et\'e 
$Y^{res}$ est \'egale au produit
fibr\'e $Y \times_X X^{res}$.

\medskip

Une r\'esolution canonique des singularit\'es 
a les propri\'et\'es suivantes d\'ecoulant de la d\'efinition~: 
\begin{enumerate}
\item [$\bullet$] les centres des \'eclatements sont 
au-dessus du lieu singulier de $X$, 
\item [$\bullet$] toute famille d'automorphismes $(\theta _g)_{g\in G}$
d'une vari\'et\'e $X$ param\'etr\'ee par une vari\'et\'e non singuli\`ere
$G$ se rel\`eve en une famille d'automorphismes de $X^{res}$.
Ceci s'applique en particulier au cas de l'action d'un groupe alg\'ebrique.
\end{enumerate}

De telles r\'esolutions canoniques des singularit\'es
existent d'apr\`es Bierstone-Milman,
Encinas-Villamayor, Hironaka ou Villamayor. 
 
Les algorithmes connus de r\'esolutions canoniques des singularit\'es
ont la propri\'et\'e sup\-pl\'em\-en\-taire suivante :
{\em pour tout faisceau d'id\'eaux coh\'erent
$I \subset {\mathcal O}_X$ sur 
une vari\'et\'e non singuli\`ere $X$, il existe  
une suite uniquement d\'etermin\'ee d'\'eclatements
le long de centres lisses $ p : X^{can} \to X$ 
de sorte que $p^{-1}(I){\mathcal O}_{X^{can}}$ 
soit principal et pour tout
morphisme lisse $f : Y \to X$, la suite correspondante d'\'eclatements
$p' : Y^{can} \to Y$ de sorte que
$(p') ^{-1} (f^* I){\mathcal O}_{Y^{can}}$ soit principal 
est \'egale au produit fibr\'e $Y \times_{X} X^{can} \to Y$.} 
(Attention, $p :X^{can} \to X$ d\'epend de $I$ mais nous n'avons
pas inclus $I$ dans la notation pour ne pas alourdir la suite.) 

Dor\'enavant, on se fixe une r\'esolution canonique des singularit\'es.
 
La proposition suivante est le coeur de cette partie :

\begin{prop}
Avec les notations du corollaire \ref{coro toro},
il y a un diagramme commutatif de vari\'et\'es alg\'ebriques :

\centerline{ 
\xymatrix{ 
& V_{i-} ^{can} \ar[dl]_{p_{i-}} \ar[d]_{h_{i-}} 
\ar@{-->}[r]^{\varphi_i^{can}} &
V_{i+} ^{can}  \ar[dr]^{p_{i+}} \ar[d]^{h_{i+}}
\\
V_i ^{res} \ar[d]_{r_i} &
 (B_{a_i}^{tor})_-/\KK^* \ar[dl]_{(\pi_{{\mathcal A}_i})_-} 
\ar@{-->}[r]^{\varphi_i^{tor}} &
(B_{a_i}^{tor})_+/\KK^* \ar[dr]^{(\pi_{{\mathcal A}_i})_+} & 
V_{i+1} ^{res} \ar[d]^{r_{i+1}} &
\\
V_i = (B_{a_i})_-/\KK^*   \ar@{-->}[rrr]^{\varphi_i} & &  & 
(B_{a_i})_+/\KK^* = V_{i+1} 
}
}

\noindent o\`u 
\begin{enumerate}
\item [(i)] $p_{i-}$ et $p_{i+}$ sont deux suites d'\'eclatements
le long de centres lisses,
\item [(ii)] si $U_{i-}^{can}= 
h_{i-}^{-1}((U_{I_{{\mathcal A}_i}}^{tor})_-/\KK^*)$ 
et $U_{i+}^{can}= h_{i+}^{-1}((U_{I_{{\mathcal A}_i}}^{tor})_+/\KK^*) $,
alors $(U_{i-}^{can} \subset V_{i-} ^{can})$ et
$(U_{i+}^{can}\subset V_{i+} ^{can})$ sont des plongements
toro\"{\i}daux et ${\varphi_i^{can}}$ est birationnelle 
toro\"{\i}dale. 
\end{enumerate}
\end{prop}

Expliquons la construction de ce diagramme.
Il s'agit de montrer dans un premier temps 
que $(\pi_{{\mathcal A}_i})_-$ 
est l'\'eclatement d'un faisceau d'id\'eaux $(I_i)_-$
sur $V_i$. C'est plus facile si l'on suppose que la famille
admissible ${\mathcal A}_i$
est choisie v\'erifiant $$\sum_{\alpha \in {\mathcal A}_i} \alpha =0,$$
ce qu'il est toujours possible de faire en rajoutant un entier
\`a la famille admissible initialement choisie. 
Alors, le faisceau toro\"{\i}dal $I_{{\mathcal A}_i}$ sur $B_{a_i}$
est engendr\'e par des fonctions invariantes 
d'apr\`es le lemme \ref{exist torique}, donc provient 
d'un faisceau d'id\'eaux $I_i$ sur $B_{a_i}/\!/\KK^*$. Le   
faisceau d'id\'eaux $(I_i)_-$ sur $V_i$ est alors obtenu
par image inverse $V_i = (B_{a_i})_-/\KK^* \to B_{a_i}/\!/\KK^*$.
Le morphisme $p_{i-} : V_{i-} ^{can} \to  V_i ^{res}$ est la suite 
d'\'eclatements uniquement d\'etermin\'ee rendant
$r_i ^{-1} ((I_i)_-) {\mathcal O}_{V_i ^{res}}$ principal
et $h_{i-}$ est l'unique morphisme induit par
la propri\'et\'e universelle de l'\'eclatement
d'un faisceau d'id\'eaux. 
La v\'erification du point (ii) se fait \`a nouveau dans
les cartes toriques fortement \'etales 
utilis\'ees pour construire le faisceau toro\"{\i}dal $I_{{\mathcal A}_i}$,
on renvoie \`a \cite{AKMW99} pour les d\'etails.

\begin{coro}\label{facto}
Soit $\varphi : X_1 \to X_2$ une application r\'eguli\`ere
birationnelle entre deux vari\'et\'es projectives lisses.
Alors, $\varphi$ se factorise en 
une suite d'\'eclatements et de contractions de centres lisses.
\end{coro}

\noindent {\em D\'emonstration.} 
C'est une cons\'equence imm\'ediate de tout ce qui pr\'ec\`ede :
il suffit d'appliquer le th\'eor\`eme \ref{factorisation toroidale}
\`a chaque ${\varphi_i^{can}}$ et de remarquer que
puisque $X_1$ et $X_2$ sont lisses, alors
$X_1 = V_0 = V_0^{res}$ et $X_2 = V_{m+1} = V_{m+1}^{res}$.\finpreuve

\subsection{Fin de la d\'emonstration du th\'eor\`eme de factorisation}
Elle se termine exactement comme dans le cas torique,
\`a l'aide du th\'eor\`eme
de Moishezon \cite{Moi67}. 
La lev\'ee des ind\'eterminations, due \`a Hironaka 
\cite{Hir75}
dans le cadre g\'en\'eral  
se formule ainsi :

\begin{theo}\label{hiro}
Soit $\varphi : X \dra X'$ une application birationnelle
entre deux vari\'et\'es com\-pl\`e\-tes lisses.
Alors il existe une 
suite d'\'eclatements le long de sous-vari\'et\'es lisses
$\psi :\bar X \to X$ 
telle que $\varphi \circ \psi : \bar X \to X'$ soit 
une application r\'eguli\`ere birationnelle.
\end{theo}

Le lemme de Chow (voir \cite{Har77} p. 107)
est lui aussi valable en g\'en\'eral :

\begin{theo}\label{chow} 
Soit $X$ une vari\'et\'e alg\'ebrique compl\`ete.
Alors il existe une vari\'et\'e projective
lisse $\tilde X$ et une application r\'eguli\`ere birationnelle
$\varphi : \tilde X \to X$. 
\end{theo}

On d\'eduit de ces r\'esultats
le th\'eor\`eme de Moishezon \cite{Moi67} 
comme dans le cas torique :

\begin{theo}\label{moishezon}
Soit $X$ une vari\'et\'e alg\'ebrique lisse et compl\`ete.
Alors il existe une vari\'et\'e projective
lisse $\tilde X$ obtenue \`a partir de $X$
par suite d'\'eclatements le long de sous-vari\'et\'es lisses
$\psi : \tilde X \to X$.
\end{theo}

Ceci permet enfin de ramener le th\'eor\`eme de 
factorisation au corollaire \ref{facto}.

\section{Appendice : 
D\'emonstration du th\'eor\`eme de $\pi$-d\'esingularisation}
\subsection{Outils et r\'esultats interm\'ediaires principaux.}
Les notes qui suivent sont une reprise de \cite{AMR99},
avec une simplification de la d\'emonstration de la 
proposition \ref{proposition fondamentale}.

\subsubsection{Relation de $\pi$-d\'ependance, $\pi$-multiplicit\'e
et $\pi$-profil de multiplicit\'e.} 
Soit $$\eta=\langle \rho _1,\ldots,\rho _k \rangle$$ 
un c\^one simplicial de $N^+$
et pour $1\leq i \leq k$, soit $v_i$ le
g\'en\'erateur dans $N$ de l'ar\^ete engendr\'ee par
$\pi (\rho _i)$. Il existe alors un unique $w_i$ dans $\QQ$
tel que $\rho _i$ soit proportionnel \`a $(v_i,w_i)$. 

\medskip

\noindent {\bf D\'efinition.}
Si $\eta$ est $\pi$-ind\'ependant, on d\'efinit la 
{\em $\pi$-multiplicit\'e} de $\eta$, not\'ee $\pimul (\eta)$
comme \'etant la multiplicit\'e du c\^one $\pi (\eta)$
(\'egale \`a l'indice dans le r\'eseau engendr\'e par $\pi (\eta)$
du sous-groupe $\ZZ v_1 \oplus \cdots \oplus \ZZ v_k$).
Remarquons que $\eta$ est $\pi$-non-singulier si et seulement
si $\pimul (\eta)=1$. 
On d\'efinit le {\em $\pi$-profil de multiplicit\'e} de $\eta$, 
not\'e $\pimp (\eta)$, 
comme \'etant le quadruplet $\pimp (\eta) := (\pimul (\eta),0,0,0)$.

\medskip

Si $\eta$ est $\pi$-d\'ependant,
comme le noyau de $\pi$ est de dimension $1$,
il existe une unique {\em relation de $\pi$-d\'ependance}
$$ \sum _{i=1}^k r_i v_i =0 \, \mbox{ avec }\,
\max \{ |r_i| \, ; \, 1\leq i \leq k \} =1  \, \mbox{ et } \, 
\sum _{i=1}^k r_i w_i > 0.$$
Pour $1\leq i \leq k$, on note $\eta _i $
la face de codimension $1$ de $\eta$ suivante~:
$$ \eta _i  = \langle 
\rho _1, \ldots, \check{\rho _i} , \ldots, \rho _k\rangle.$$
 Le lemme suivant est \'el\'ementaire et essentiel~:

\begin{lemm}   
\begin{enumerate}
\item [(i)] La face $\eta _i$ est $\pi$-ind\'ependante si et seulement 
si $r_i \neq 0$.
\item [(ii)] Si $\eta _i$ et $\eta _j$ sont toutes deux 
$\pi$-ind\'ependantes, alors 
$|r_j| \pimul (\eta _i) = |r_i| \pimul (\eta _j).$
\end{enumerate}
\end{lemm}

\noindent {\bf Notation.}
Introduisons les ensembles~: 
\noindent $I_1(\eta)  = \{ i \, ; \, r_i =1\} \, ,\,  
I_+(\eta)  = \{ i \, ; \, r_i >0\}\,,\,
I_{-1}(\eta)  = \{ i \, ; \, r_i = -1\}\,,\,
I_{-}(\eta)  = \{ i \, ; \, r_i < 0\}$
et $i_1(\eta)$, $i_+(\eta)$, $i_{-1}(\eta)$ et $i_-(\eta)$ 
leur cardinal respectif.

\medskip

\noindent {\bf D\'efinition.}
Si $\eta$ est $\pi$-d\'ependant, on d\'efinit la 
{\em $\pi$-multiplicit\'e} de $\eta$
comme \'etant le maximum des $\pi$-multiplicit\'es des
faces $\eta _i$ pour $i \in I_+(\eta) \cup I_{-}(\eta)$.
Autrement dit, $\pimul (\eta)$ est le maximum des $\pi$-multiplicit\'es des
faces $\pi$-ind\'ependantes de $\eta$. 
On d\'efinit le {\em $\pi$-profil de multiplicit\'e} de $\eta$, 
not\'e $\pimp (\eta)$, 
comme \'etant le quadruplet 
$$\pimp (\eta) := (\pimul (\eta),0,0,0)\, \mbox { si }\,
i_1(\eta) + i_{-1}(\eta) =1,$$ et 
$$\pimp (\eta) := (\pimul (\eta),1,i_+(\eta) + i_{-}(\eta),
i_1(\eta) + i_{-1}(\eta))\, \mbox { si }\,
i_1(\eta) + i_{-1}(\eta) \geq 2.$$

\medskip

Remarquons que la 
quantit\'e $k_{\eta} = i_+(\eta) + i_{-}(\eta)$ est la dimension
de l'unique circuit contenu dans $\eta$.

\medskip

\noindent {\bf Notations.} 
Soit $\Sigma$ un \'eventail de $N^+_{\QQ}$.
On note $g_{\Sigma}$ le maximum (pour l'ordre lexicographique)
des $\pimp (\eta)$ lorsque $\eta$
d\'ecrit tous les c\^ones maximaux de $\Sigma$ ({\em i.e.} non contenus
strictement dans un c\^one de $\Sigma$)
et $s_{\Sigma}$ le nombre de c\^ones maximaux atteignant ce maximum.
Alors le {\em $\pi$-profil
de multiplicit\'e} de $\Sigma$ est le quintuplet~:
$$ \pimp (\Sigma) := (g_{\Sigma},s_{\Sigma}).$$

\subsubsection{Faces cod\'efinies.}
Soit $\eta = 
\langle\rho _1,\ldots,\rho _k\rangle$ un c\^one 
$\pi$-d\'ependant de $N^+_{\QQ}$
et $\tau$ une face $\pi$-ind\'ependante de $\eta$.
Si
$$ \sum _{i=1}^k r_i v_i =0 $$ est la relation de $\pi$-d\'ependance
de $\eta$,
on dit que $\tau$ est {\em cod\'efinie} par rapport
\`a $\eta$ si ses g\'en\'erateurs sont inclus dans 
$\{\rho _i \,,\, r_i \geq 0\}$
ou dans $\{\rho _i \,,\, r_i \leq 0\}$.
Evidemment, si $\tau$ est une face de $\tau '$ 
o\`u $\tau '$ est cod\'efinie par rapport
\`a $\eta$, alors $\tau$ est cod\'efinie par rapport
\`a $\eta$.

Si $\tau = \langle\rho _1,\ldots,\rho _l\rangle$ 
est un c\^one $\pi$-ind\'ependant de $N^+_{\QQ}$,
et si  $v_i$ est le
g\'en\'erateur dans $N$ de l'ar\^ete engendr\'ee par
$\pi (\rho _i)$, on note~: 
$$ \para (\pi (\tau)) = 
\{ v \in N \, ; \, v = \sum _{i=1}^l a_iv_i \,, \, 0< a_i < 1 \}.$$
Remarquons que si $\eta$ est un c\^one $\pi$-singulier de $N^+_{\QQ}$,
et si $\tau$ est une face $\pi$-ind\'ependante et $\pi$-singuli\`ere 
de $\eta$ de dimension
minimale, alors $\para (\pi (\tau)) \neq \emptyset$.

La proposition suivante est le r\'esultat central de
ce paragraphe. Elle
illustre parfaitement l'utilit\'e des faces cod\'efinies.

\begin{prop}\label{codefinie}
Soit $\eta$ un c\^one $\pi$-d\'ependant de $N^+$
et $\tau$ une face $\pi$-ind\'ependante 
de $\eta$. Soit $v \in \para (\pi (\tau))$
et $\rho \in N_{\QQ}^+$ dans l'int\'erieur relatif de 
$\tau$ tel que $\pi (\rho) = v$. Soit $\eta '$
l'\'eventail obtenu par subdivision \'etoil\'ee de $\eta$
par rapport \`a $\rho$.
Si $\tau$ est cod\'efinie par rapport \`a
$\eta$, alors $\pimp (\eta ') \leq \pimp (\eta)$
et si de plus $\tau$ est contenue dans une 
face $\gamma$ de codimension $1$ de $\eta$ de $\pi$-multiplicit\'e
maximale, alors $\pimp (\eta ') < \pimp (\eta)$.
\end{prop}

\noindent {\em D\'emonstration.} Notons
$\tau = \langle\rho_1,\ldots,\rho_m\rangle$ et 
$\eta = \langle\rho_1,\ldots,\rho_n\rangle$ avec $n>m$.
On note $\rho = \sum_{i=1}^m a_i\rho_i$
et soit $\sum_{i=1}^n r_i v_i=0$ la relation de $\pi$-d\'ependance de
$\eta$. 

Les c\^ones maximaux de $\eta '$ sont les
$\eta _{\alpha} = 
\langle\rho,\rho_1, \ldots ,\check{ \rho_{\alpha}},\ldots,\rho_m,\ldots,\rho_n\rangle$
pour $1\leq \alpha \leq m$.
Nous allons estimer $\pimp(\eta _{\alpha})$.

Les faces nouvelles de $\eta _{\alpha}$ ({\em i.e.} 
les faces de $\eta _{\alpha}$ 
qui ne sont pas face de $\eta$)
sont les
$$\gamma_{\alpha \beta} = 
\langle\rho,\rho_1,\ldots,\check{\rho_{\alpha}},\ldots,
\check{\rho_{\beta}},\ldots,\rho_n\rangle$$
pour $\beta \neq \alpha$. Supposons que $\gamma_{\alpha \beta}$
est $\pi$-ind\'ependante et remarquons que cela implique que
$r_{\alpha}\neq 0$ ou $r_{\beta}\neq 0$.  
\begin{enumerate}
\item[$\bullet$] si $\beta \geq m+1$, alors 
$$
\begin{array}{rcl}
\pimul (\gamma_{\alpha \beta})&=&|\det (v,v_1,
\ldots,\check{v_{\alpha}},\ldots,
\check{v_{\beta}},\ldots,v_n) |\\
& = & 
a_{\alpha}|\det (v_1,
\ldots,\check{v_{\beta}},\ldots,v_n) | \\
& \leq & a_{\alpha}\pimul (\eta) <
\pimul (\eta).
\end{array}
$$ 
\item[$\bullet$] si $\beta \leq m$, alors
$$
\begin{array}{rcl}
\pimul (\gamma_{\alpha \beta})&=&
|a_{\alpha}\det (v_{\alpha},v_1,\ldots,\check{v_{\alpha}},
\ldots,\check{v_{\beta}},\ldots,v_n) \\
& & +a_{\beta}\det (v_{\beta},v_1,\ldots,\check{v_{\alpha}},
\ldots,\check{v_{\beta}},\ldots,v_n)|.                        
\end{array}
$$
Pour $\lambda = \alpha$ ou $\beta$, posons 
$m_{\lambda} = \det (v_{\lambda},v_1,\ldots,\check{v_{\alpha}},
\ldots,\check{v_{\beta}},\ldots,v_n)$.
De $\sum_{i=1}^n r_i v_i=0$, on en d\'eduit 
que $r_{\alpha}m_{\alpha} +r_{\beta}m_{\beta} =0$, 
et $\tau$ \'etant cod\'efinie par rapport \`a $\eta$,
$r_{\alpha}$ et $r_{\beta}$ sont de m\^eme signe
donc $m_{\alpha}$ et $m_{\beta}$ sont de signes oppos\'es.
De l\`a
$$
\begin{array}{rcl} 
\pimul (\gamma_{\alpha \beta})&=&| a_{\alpha} m_{\alpha}+a_{\beta} 
m_{\beta}|\\
&\leq& \max (a_{\alpha} |m_{\alpha}|,a_{\beta} |m_{\beta}|) \\
&\leq& \max (a_{\alpha},a_{\beta}) \pimul (\eta) < \pimul (\eta).
\end{array}
$$

\noindent Ainsi, pour toute nouvelle face $\gamma_{\alpha \beta}$,
on a $\pimul (\gamma_{\alpha \beta}) < \pimul (\eta)$.
\end{enumerate}

Enfin, chaque $\eta _{\alpha}$ poss\`ede une unique ancienne face
de codimension un~: $$\gamma _{\alpha} = 
\langle\rho_1,\ldots,\check{\rho_{\alpha}},\ldots,\rho_n\rangle$$ et cette 
face satisfait \'evidemment
$\pimul (\gamma _{\alpha})\leq \pimul (\eta)$.

Le bilan est le suivant~: si 
$\pimul (\gamma _{\alpha}) =  \pimul (\eta)$,
alors $\pimp (\eta_{\alpha}) = (\pimul (\eta),0,0,0)$
et si $\pimul (\gamma _{\alpha}) <  \pimul (\eta)$, alors
$\pimp (\eta_{\alpha}) < (\pimul (\eta),0,0,0)$. 

Notons $s'$ le nombre de faces de codimension $1$
de $\eta$ ne contenant pas $\tau$ et de $\pi$-multiplicit\'e maximale
et $r_{\eta}$ le nombre de faces de codimension $1$
de $\eta$ de $\pi$-multiplicit\'e maximale.

Alors, 
\begin{enumerate}
\item[$\bullet$]
si $s'=0$, ce qui pr\'ec\`ede montre que 
$\pimp (\eta_{\alpha}) < (\pimul (\eta),0,0,0)$
pour tout $\alpha$, $1\leq \alpha \leq m$, 
et donc $\pimp (\eta ') < \pimp (\eta)$.
\item[$\bullet$] si $s'=1$, on a 
$$\pimp (\eta ') = ((\pimul (\eta),0,0,0),1) \leq \pimp (\eta).$$
De plus, cette in\'egalit\'e est une \'egalit\'e si et seulement
si $r_{\eta}=s'=1$, ce qui ne se produit pas dans le cas o\`u
$\tau$ est contenue dans une face de codimension $1$
de $\eta$ de $\pi$-multiplicit\'e maximale.
\item[$\bullet$] 
si $s'\geq 2$, on a 
$$\pimp (\eta ') = ((\pimul (\eta),0,0,0),s') 
< \pimp (\eta) = ((\pimul (\eta),1,k_{\eta},r_{\eta}),1).$$
\end{enumerate}
\finpreuve

\medskip

La proposition \ref{codefinie} nous encourage \`a fabriquer des faces
cod\'efinies. Ce sera l'objet du paragraphe suivant.
Enon\c cons d\`es maintenant le lemme imm\'ediat suivant~:

\begin{lemm}\label{lemme facile}
Soient $\eta$ un c\^one $\pi$-d\'ependant de $N^+$,
$\tau$ une face $\pi$-ind\'ependante de $\eta$ et $\sigma$
l'unique circuit contenu dans $\eta$.
Si $\dim (\sigma) \leq 2$, alors 
$\tau$ est cod\'efinie par rapport \`a
$\eta$.
\end{lemm}

\subsubsection{Subdivision \'etoil\'ee positive et n\'egative.}
Soient $\Sigma$ un \'eventail de $N_{\QQ}^+$
et $\sigma = \langle\rho _1,\ldots,\rho _k\rangle$ 
un circuit de $\Sigma$. 
Si
$$ \sum _{i=1}^k r_i v_i =0 $$ est la relation de $\pi$-d\'ependance
de $\sigma$, on pose
$$ v_+ = \sum_{i \in I_+(\sigma)} v_i \, \mbox{ et } \,
v_- = \sum_{i \in I_-(\sigma)} v_i .$$
(Attention, $v_+$ et $v_-$ ne sont pas n\'ecessairement primitifs.)
Remarquons que comme $\sigma$ est un circuit,
on a $I_+(\sigma) \cup I_-(\sigma) = \{1,\ldots,k\}$, 
autrement dit, tous les $r_i$ de la relation de $\pi$-d\'ependance
de $\sigma$ sont non nuls.

Comme $$ v_+ = \sum_{i \in I_+(\sigma)} v_i = 
\sum_{i \in I_+(\sigma)} (1 - \varepsilon r_i)v_i + 
\sum_{i \in I_-(\sigma)} (- \varepsilon r_i)v_i,$$
pour $\varepsilon$ petit, $ v_+$
est de la forme $\sum_{i=1}^k c_i v_i$ o\`u
tous les $c_i$ sont strictement positifs.
Il existe donc $\rho _+$ dans l'int\'erieur relatif de 
$\sigma$ tel que $\pi (\rho _+) = v_+$. Soit $\Sigma^+$
l'\'eventail obtenu par subdivision \'etoil\'ee de $\Sigma$
par rapport \`a $\rho _+$. On l'appelle {\em subdivision 
\'etoil\'ee positive} de $\Sigma$ par rapport \`a $\sigma$.
(Cette construction d\'epend du choix de $\rho _+$
mais les propri\'et\'es que nous \'enoncerons n'en d\'ependent
pas.) La m\^eme construction \`a partir de $v_-$
donne lieu \`a la {\em subdivision 
\'etoil\'ee n\'egative} $\Sigma^-$ de $\Sigma$ par rapport \`a $\sigma$.

Malheureusement, le $\pi$-profil de multiplicit\'e
de $\Sigma^+$ (resp. $\Sigma^-$) n'est en g\'en\'eral pas inf\'erieur
ou \'egal \`a celui de $\Sigma$.
Le lemme facile suivant illustre n\'eanmoins l'int\'er\^et 
des subdivisions \'etoil\'ees positive et n\'egative~:

\begin{lemm}\label{lemme codefini}
Soient $\eta$ un c\^one $\pi$-d\'ependant de $N^+_{\QQ}$,
$\gamma$ une face de codimension $1$ de $\eta$ de $\pi$-multiplicit\'e
maximale. 
Soit $\eta^+$ (resp. $\eta^-$)
l'\'eventail obtenu par subdivision \'etoil\'ee positive 
(resp. n\'egative) de $\eta$
par rapport \`a  l'unique circuit $\sigma$ contenu dans $\eta$. 
Alors $\gamma$ est cod\'efinie par rapport au
c\^one maximal de $\eta^+$ (resp. $\eta^-$) contenant $\gamma$. 
\end{lemm} 

\medskip

La proposition suivante joue
un r\^ole crucial dans la d\'emonstration du th\'eor\`eme
de $\pi$-d\'esin\-gula\-risa\-tion.

\begin{prop}\label{proposition fondamentale}
Soit $\sigma$
un circuit de $N^+_{\QQ}$ de dimension strictement sup\'erieure \`a $2$.
Alors soit l'\'eventail $\sigma ^+$ obtenu par subdivision 
\'etoil\'ee positive 
de $\sigma$, soit l'\'eventail $\sigma ^-$
obtenu par subdivision \'etoil\'ee n\'egative 
de $\sigma$ satisfait l'une des deux propri\'et\'es suivantes
(nous notons $\sigma '$ l'\'eventail $\sigma ^+$ ou 
l'\'eventail $\sigma ^-$ qui convient)~:
\begin{enumerate}
\item [A)] pour tout c\^one maximal $\delta '$
de $\sigma '$, on a $\pimp (\delta ') < \pimp (\sigma)$~;
en particulier $$\pimp (\sigma ') < \pimp (\sigma).$$
\item [B)] il existe un c\^one maximal $\kappa '$
de $\sigma '$ v\'erifiant $\pimp (\kappa ') = \pimp (\sigma)$ et
tel que tout c\^one maximal $\delta ' \neq \kappa '$
de $\sigma '$ v\'erifie
$\pimp (\delta ') < \pimp (\sigma)$~; en particulier 
$$\pimp (\sigma ') = \pimp (\sigma).$$ 
De plus, l'unique ancienne 
face $\gamma '$ de codimension $1$ de $\kappa '$ est $\pi$-ind\'ependante,
et v\'erifie 
$\pimul (\gamma ') = \pimul (\kappa ') = \pimul (\sigma)$.
En particulier, $\gamma '$ est cod\'efinie par rapport
\`a $\kappa '$ d'apr\`es le lemme \ref{lemme codefini}. 
\end{enumerate}
\end{prop}

Avant de commencer la d\'emonstration 
proprement dite, commen\c cons par un calcul interm\'e\-diai\-re~: notons
$\sigma = \langle\rho_1,\ldots,\rho_n\rangle$, 
$\sum_{i=1}^n r_i v_i=0$ la relation de $\pi$-d\'ependance de
$\sigma$ et $v_+ = \sum_{i \in I_+(\sigma)} v_i = e\bar{v}_+$
avec $\bar{v}_+$ primitif. Soit enfin   
$\rho _+$ dans l'int\'erieur relatif de 
$\sigma$ tel que $\pi (\rho _+) = v_+$.
Les c\^ones maximaux de $\sigma ^+$
sont de la forme $\sigma _{\alpha} = 
\langle\rho _+,\rho_1,\ldots,\check{\rho_{\alpha}},\ldots,\rho_n\rangle$
et nous devons estimer leur $\pi$-profil de multiplicit\'e~;
chacun d'eux a une face ancienne de codimension un
$\gamma_{\alpha} = \langle\rho_1,\ldots,\check{\rho_{\alpha}},\ldots,
\rho_n\rangle$
et de nouvelles faces 
$$\gamma_{\alpha \beta} = \langle\rho_+,\rho_1,\ldots,\check{\rho_{\alpha}},\ldots,
\check{\rho_{\beta}},\ldots,\rho_n\rangle.$$
Les calculs de $\pi$-multiplicit\'e des faces
de $\sigma _{\alpha}$ sont faciles :

\begin{enumerate}
\item[$\bullet$] si $r_{\alpha} > 0$ alors 

$$
\pimul (\gamma_{\alpha \beta}) = 
\left\{
\begin{array}{l} \frac{1}{e}|\pimul (\gamma _{\alpha})-
\pimul (\gamma _{\beta})| \, \, \mbox{ si } \, r_{\beta} > 0 \\
\\
\frac{1}{e} \pimul (\gamma _{\beta}) \, \, \mbox{ si } \, r_{\beta} < 0,\\ 
\end{array}
\right.
$$
\item[$\bullet$] si $r_{\alpha} < 0$ alors 

$$
\pimul (\gamma_{\alpha \beta}) = 
\left\{
\begin{array}{l} \frac{1}{e}\pimul (\gamma _{\alpha})
\, \, \mbox{ si } \, r_{\beta} > 0 \\
\\
0  \, \, \mbox{ si } \, r_{\beta} < 0.\\ 
\end{array}
\right.
$$
\end{enumerate}

\medskip

\noindent {\em D\'emonstration de la proposition \ref{proposition
fondamentale}.}

\begin{enumerate}
\item[$\bullet$] supposons qu'une seule face de codimension~1
de $\sigma$ est
de $\pi$-multiplicit\'e maximale, autrement dit
$\pimp (\sigma) = (\pimul (\sigma),0,0,0)$.
Ceci signifie que $i_1(\sigma) + i_{-1}(\sigma)=1$
et on peut supposer que $i_{-1}(\sigma)=0$ et $i_1(\sigma)=1$, on
note alors $\alpha_0$ le seul indice tel que $r_{\alpha_0}=1$. 
Montrons que $\sigma^+$
satisfait la propri\'et\'e B). 

En effet, les formules pr\'ec\'edentes montrent que  
$\pimul (\sigma _{\alpha}) < \pimul (\sigma)$ pour $\alpha \neq \alpha _0$
(et donc $\pimp (\sigma _{\alpha}) < 
\pimp (\sigma)$ pour $\alpha \neq \alpha _0$)
et que $\pimul (\sigma _{\alpha _0}) = \pimul (\sigma)$, la seule face
de codimension $1$ de $\sigma _{\alpha _0}$ de 
$\pi$-multiplicit\'e maximale \'etant 
$\gamma _{\alpha _0}$. On en d\'eduit que $\pimp (\sigma _{\alpha _0})
= \pimp (\sigma)$, et $\gamma _{\alpha _0}$, \'etant de 
$\pi$-multiplicit\'e maximale, est cod\'efinie par rapport \`a 
$\sigma _{\alpha _0}$
d'apr\`es le lemme \ref{lemme codefini}. 

\item[$\bullet$] supposons que 
$\pimp (\sigma) = 
(\pimul (\sigma),1,i_+(\sigma)+i_-(\sigma),i_1(\sigma) + i_{-1}(\sigma))$.
Dans ce cas, $i_1(\sigma) + i_{-1}(\sigma) \geq 2$ et on peut
supposer que $i_1(\sigma) \geq 1$. 
Montrons que $\sigma^+$
satisfait l'une des propri\'et\'es A) ou B) sauf si
$i_{-1}(\sigma)=i_{-}(\sigma)=1$.  

\begin{enumerate}
\item[(i)] si $e>1$, les formules pr\'ec\'edentes montrent que
$\pimul (\sigma _{\alpha}) \leq \pimul (\sigma)$ pour tout $\alpha$,
et que si $\pimul (\sigma _{\alpha}) = \pimul (\sigma)$,
la seule face
de codimension $1$ de $\sigma _{\alpha}$ de 
$\pi$-multiplicit\'e maximale est 
$\gamma _{\alpha}$. On en d\'eduit que 
$\pimp (\sigma _{\alpha}) \leq (\pimul(\sigma),0,0,0) < 
\pimp (\sigma)$ pour tout $\alpha$
et donc que $\sigma^+$
satisfait la propri\'et\'e A).     

\item[(ii)] si $e=1$, \'evaluons s\'epar\'ement $\pimp (\sigma _{\alpha})$
suivant que $r_{\alpha}=-1$,
$-1< r_{\alpha} <0$, $r_{\alpha} =1$ ou $1> r_{\alpha} >0$.
\begin{enumerate}
\item[a)] si $-1< r_{\alpha} <0$, 
les formules pr\'ec\'edentes montrent que
$\pimul (\sigma_{\alpha}) < \pimul (\sigma)$ et donc 
que $ \pimp (\sigma_{\alpha})< \pimp (\sigma)$.
\item[b)] si $1> r_{\alpha} >0$, les formules pr\'ec\'edentes montrent que 
$$
\pimp (\sigma_{\alpha}) 
\left\{
\begin{array}{cl} \leq &(\pimul(\sigma),0,0,0) < \pimp(\sigma)
\, \, \mbox{ si } \, i_{-1}({\sigma})\leq 1 ,\\
\\
= & (\pimul(\sigma),1,k_{\sigma_{\alpha}}, i_{-1}({\sigma}))
\, \, \mbox{ si } \, i_{-1}({\sigma})\geq 2. 
\end{array}
\right.
$$
Dans le deuxi\`eme cas, $k_{\sigma_{\alpha}}\leq k_{\sigma}=\dim(\sigma)$
et $i_{-1}({\sigma}) < i_1(\sigma) + i_{-1}(\sigma)$
donc $\pimp ( \sigma_{\alpha}) < \pimp ( \sigma)$.

\item[c)] si $r_{\alpha}=1$, les formules pr\'ec\'edentes montrent que
$$
\pimp (\sigma_{\alpha}) = 
\left\{
\begin{array}{l} (\pimul(\sigma),1,k_{\sigma_{\alpha}}, 
1+i_{-1}({\sigma}))
\, \, \mbox{ si } \, i_{-1}({\sigma})\neq 0,\\
(\pimul(\sigma),0,0,0) < \pimp(\sigma)
\, \, \mbox{ si } \, i_{-1}({\sigma}) =0 .\\
\end{array}
\right.
$$
Dans le premier cas, 
si $i_{1}(\sigma)\geq 2$, $\sigma_{\alpha}$ poss\`ede une face 
de codimension $1$ $\pi$-d\'ependante 
(car de $\pi$-multiplicit\'e nulle)
donc $k_{\sigma_{\alpha}} < k_{\sigma}$ et par suite
$\pimp (\sigma_{\alpha}) < \pimp (\sigma)$~;
sinon $i_{1}(\sigma)=1$, par suite  
$\pimp (\sigma_{\alpha}) = \pimp (\sigma)$ et 
$\gamma_{\alpha}$, \'etant de 
$\pi$-multiplicit\'e maximale, est cod\'efinie par rapport \`a 
$\sigma _{\alpha}$
d'apr\`es le lemme \ref{lemme codefini}.

\item[d)] si $r_{\alpha}=-1$, les formules pr\'ec\'edentes montrent que
$$
\pimp (\sigma_{\alpha}) = 
(\pimul(\sigma),1,k_{\sigma_{\alpha}}, 
1+i_{+}({\sigma})).$$ 
(Remarquons que $i_{+}({\sigma})\geq i_{1}({\sigma})\geq 1$).
Si $i_{-}({\sigma}) \geq 2$, 
$\sigma_{\alpha}$ poss\`ede une face 
de codimension $1$ $\pi$-d\'ependante (car de $\pi$-multiplicit\'e nulle)
donc $k_{\sigma_{\alpha}} < k_{\sigma}$ et par suite
$\pimp (\sigma_{\alpha}) < \pimp (\sigma)$~; 
si $i_{-}({\sigma}) =1 $, alors
$$ 1+i_{+}({\sigma}) \geq i_1(\sigma) + i_{-1}(\sigma) = i_1(\sigma)+1 $$
et on ne peut conclure \`a ce stade de la d\'emonstration. 
\end{enumerate}

\smallskip

Le bilan est cependant le suivant~: si 
$i_1(\sigma) + i_{-1}(\sigma) \geq 2$ et $i_1(\sigma) \geq 1$
alors $\sigma^+$
satisfait l'une des propri\'et\'es A) ou B)
sauf si $i_{-1}(\sigma)=i_{-}(\sigma)=1$.
De fa\c con sym\'etrique, si 
$i_1(\sigma) + i_{-1}(\sigma) \geq 2$ et $i_{-1}(\sigma) \geq 1$
alors $\sigma^-$
satisfait l'une des propri\'et\'es A) ou B)
sauf si $i_{1}(\sigma)=i_{+}(\sigma)=1$.

Le seul cas restant est donc celui o\`u 
$$i_{-1}(\sigma)=i_{-}(\sigma)=i_{1}(\sigma)=i_{+}(\sigma)= 1,$$
mais alors $\dim(\sigma) =2$ ce qui est exclu par l'hypoth\`ese.
\end{enumerate}
\end{enumerate}
\finpreuve

\medskip

Nous aurons besoin du lemme suivant qui 
pr\'ecise la proposition \ref{proposition fondamentale} :

\begin{lemm}\label{lemme fondamental}
Soit $\Sigma$ un \'eventail de $N^+$, 
$\sigma$
un circuit de $\Sigma$ de dimension strictement sup\'erieure \`a $2$.
Si $\overline{\st (\sigma)}'$ d\'esigne
l'\'eventail obtenu par la subdivision \'etoil\'ee positive ou n\'egative
de $\overline{\st (\sigma)}$ donn\'ee par la proposition pr\'ec\'edente, 
alors $$\pimp (\overline{\st (\sigma)}') <
\pimp (\overline{\st (\sigma)})\, \mbox { dans le cas A) } $$
et 
$$\pimp (\overline{\st (\sigma)}') =
\pimp (\overline{\st (\sigma)})\, \mbox { dans le cas B) } .$$
\end{lemm}

\noindent {\em Indication de d\'emonstration.} 
On montre que si 
$$\pimp (\sigma) = (\pimul(\sigma),b_{\sigma},k_{\sigma},
r_{\sigma},1) \, \mbox{ et } \,
\pimp (\sigma ') = (\pimul(\sigma '),b_{\sigma '},k_{\sigma '},
r_{\sigma '},s'),$$ alors il existe des entiers $e$ et $s$
sup\'erieurs ou \'egaux \`a $1$ tels que
$$\pimp (\overline{\st (\sigma)}) = (\pimul(\sigma)e,b_{\sigma},k_{\sigma},
r_{\sigma},s)$$ et 
$$\pimp (\overline{\st (\sigma)}') = (\pimul(\sigma ')e,b_{\sigma '},
k_{\sigma '},r_{\sigma '},s's).$$  
\finpreuve

\subsection{D\'emonstration du th\'eor\`eme de $\pi$-d\'esingularisation.}
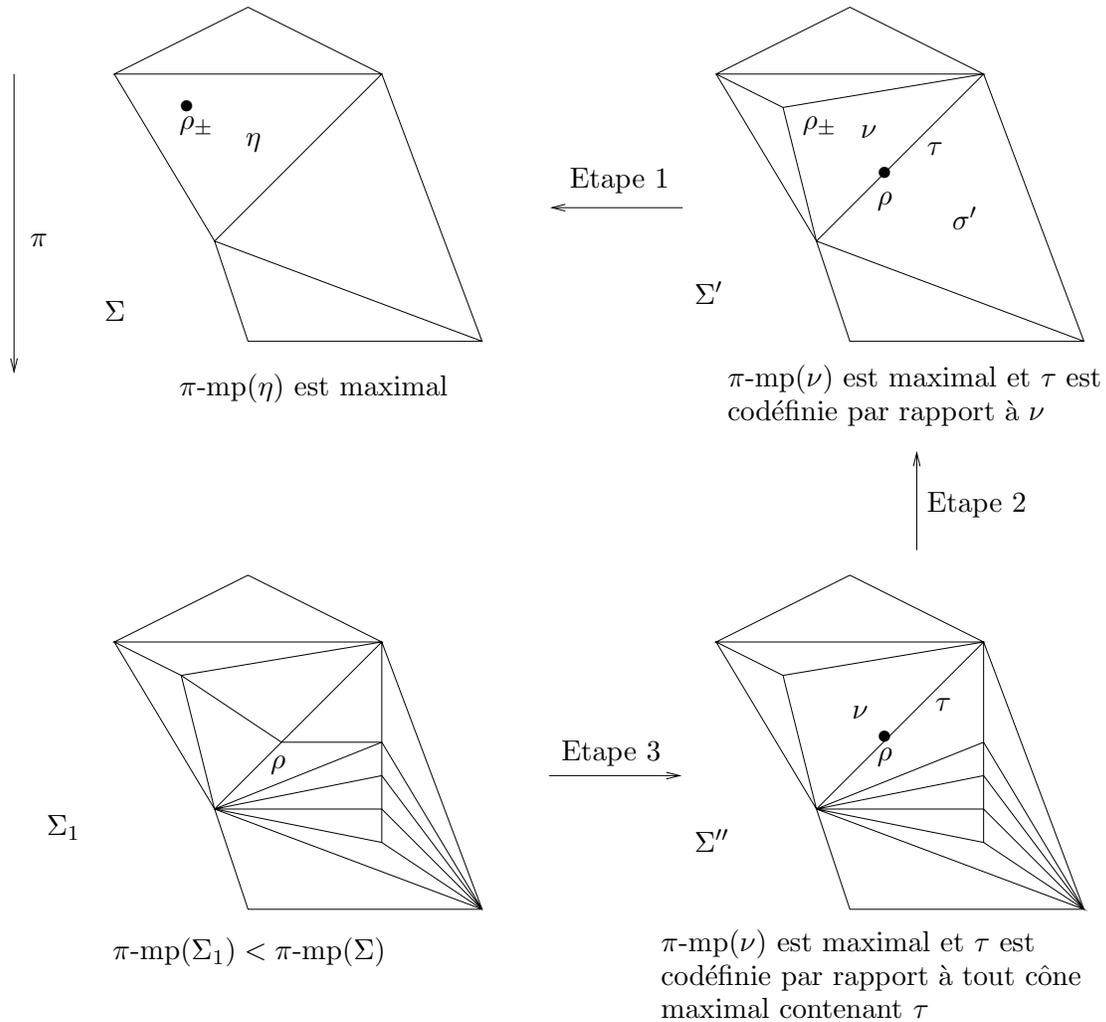
\begin{figure}[ht]
  \begin{center}
    \leavevmode
    \input{fig2.pstex_t}
    \caption{Le th\'eor\`eme de $\pi$-d\'esingularisation}
    \label{figure fondamentale}
  \end{center}
\end{figure}

Soit $\Sigma$ un \'eventail de $N^+$.
La strat\'egie est claire~: si $\Sigma$ est $\pi$-non-singulier,
il n'y a rien \`a faire. Sinon, il suffit
de construire un \'eventail $\Sigma_1$ 
obtenu par une suite 
finie de subdivisions
\'etoil\'ees de $\Sigma$ n'affectant pas 
les c\^ones $\pi$-non-singuliers de $\Sigma$
et v\'erifiant $\pimp(\Sigma _1) < \pimp(\Sigma)$.
Si $\Sigma_1$ est $\pi$-non-singulier, c'est fini, sinon
on recommence. Ce proc\'ed\'e doit s'arr\^eter apr\`es un nombre fini
d'\'etapes et l'\'eventail obtenu est $\pi$-non-singulier. 
La construction de $\Sigma_1$ se fait en trois \'etapes, illustr\'ees
par la figure \ref{figure fondamentale}.

{\bf Etape 1.}
Posons $\pimp(\Sigma)=
(g_{\Sigma}, s_{\Sigma})$ et choisissons
$\eta$ un c\^one maximal de $\Sigma$ de $\pi$-profil de multiplicit\'e
maximal ({\em i.e.} \'egal \`a $g_{\Sigma}$) et soit $\sigma$
l'unique circuit contenu dans $\eta$.

\begin{enumerate}
\item[(i)]
Supposons que $\dim (\sigma) > 2$ et
appliquons alors la proposition \ref{proposition
fondamentale}~; on note $\Sigma '$
l'\'eventail ainsi obtenu et
$\eta '$ le sous-\'eventail de $\Sigma '$
obtenu en subdivisant $\eta$. 
\begin{enumerate}
\item[$\bullet$]
Dans le cas A),
le lemme \ref{lemme fondamental} 
assure que $\pimp (\Sigma ') < \pimp (\Sigma)$.
On pose alors $\Sigma_1 := \Sigma$.
\item[$\bullet$]
Dans le cas B), le lemme \ref{lemme fondamental} assure que 
$\pimp (\Sigma ') = \pimp (\Sigma)$. Soit $\gamma$ l'unique
face de codimension $1$ de $\eta$ telle que 
$\gamma \cap \sigma = \gamma '$ ($\gamma '$ donn\'ee par
la proposition \ref{proposition
fondamentale}). Alors $\gamma$
est cod\'efinie par rapport \`a l'unique c\^one maximal 
$\nu$
de $\eta '$ dont elle est face, elle est de plus 
de $\pi$-multiplicit\'e maximale.
\end{enumerate}

\item[(ii)] Supposons que $\dim (\sigma) \leq 2$ et
choisissons une 
face $\gamma$ de codimension $1$ de $\eta$,
de $\pi$-multiplicit\'e maximale. Comme $\dim (\sigma) \leq 2$,
$\gamma$ est cod\'efinie par rapport \`a $\eta$.
\end{enumerate}

Le bilan de l'Etape~1 est le suivant~: nous avons
construit un \'eventail $\Sigma '$, subdivision \'etoil\'ee de $\Sigma$,
v\'erifiant
$\pimp (\Sigma ') \leq \pimp (\Sigma)$ avec 
\begin{enumerate}
\item[(i)] un
c\^one maximal $\nu$ de $\pi$-profil de multiplicit\'e maximal,
\item[(ii)] une face $\gamma$ de codimension $1$ de $\nu$,
de $\pi$-multiplicit\'e maximale, cod\'efinie par rapport \`a~$\nu$. 
\end{enumerate}

Choisissons alors  
$\tau$ une face de $\gamma$ (donc cod\'efinie par rapport
\`a $\nu$), $\pi$-singuli\`ere
de dimension minimale, $v \in \para (\pi (\tau))$
et $\rho \in N_{\QQ}^+$ dans l'int\'erieur relatif de 
$\tau$ tel que $\pi (\rho) = v$.

Le probl\`eme \`a ce stade est que $\tau$, qui est cod\'efinie par rapport
\`a $\nu$, ne l'est en g\'en\'eral pas par rapport aux autres 
c\^ones maximaux qui la contiennent.

{\bf Etape 2.}
Enon\c cons la proposition suivante, dont la d\'emonstration sera 
donn\'ee \`a la fin de ce paragraphe~:

\begin{prop}\label{prop finale}
Soit $\Sigma$ un \'eventail de $N^+$, $\sigma$
un circuit de $\Sigma$.
Soit $\tau$ dans $\overline{\st (\sigma)}$, 
$\pi$-ind\'ependante.
Alors il y a une subdivision de $\overline{\st (\sigma)}$,
not\'ee $\overline{\st (\sigma)}'$, obtenue par une suite
finie de subdivisions \'etoil\'ees positives ou n\'egatives 
par rapport \`a des circuits successifs
contenus dans $\sigma$ telle que~:
\begin{enumerate}
\item[(i)] $\pimp (\overline{\st (\sigma)}') \leq 
\pimp (\overline{\st (\sigma)}),$
\item[(ii)] $\tau$ est une face de $\overline{\st (\sigma)}'$ 
({\em i.e.} $\tau$ n'est pas affect\'ee par les subdivisions)
et $\tau$ est cod\'efinie par rapport \`a tous les c\^ones 
maximaux de $\overline{\st (\sigma)}'$ qui la contiennent.
\end{enumerate}
\end{prop}

Appliquons cette proposition de la fa\c con suivante
\`a notre situation~: 
on consid\`ere l'ensemble des circuits $\sigma '$
diff\'erents de celui contenu dans $\nu$ 
tels que $\tau$ est contenu dans $\overline{\st (\sigma ')}$
et on applique la proposition pr\'ec\'edente 
\`a chacun d'eux (le r\'esultat obtenu ne d\'epend pas 
de l'ordre dans lequel on a consid\'er\'e les diff\'erents $\sigma'$). 
 
\smallskip

Le bilan de l'Etape~2 est alors le suivant~: nous avons
construit un \'eventail $\Sigma ''$, 
obtenu par une suite de subdivisions \'etoil\'ees de $\Sigma$
v\'erifiant
$\pimp (\Sigma '') \leq \pimp (\Sigma)$, avec 
\begin{enumerate}
\item[(i)]
un c\^one maximal $\nu$ de $\pi$-profil de multiplicit\'e maximal,
\item[(ii)] une face $\gamma$ de codimension $1$ de $\nu$,
de $\pi$-multiplicit\'e maximale, cod\'efinie par rapport \`a~$\nu$,  
\item[(iii)] une face
$\tau$ de $\gamma$ (donc cod\'efinie par rapport
\`a~$\nu$), $\pi$-singuli\`ere
de dimension minimale, $v \in \para (\pi (\tau))$
et $\rho \in N_{\QQ}^+$ dans l'int\'erieur relatif de 
$\tau$ tel que $\pi (\rho) = v$.
{\em  De plus, $\tau$ est cod\'efinie par rapport
\`a tous les c\^ones 
maximaux de $\Sigma ''$ qui la contiennent.}   

\end{enumerate}

{\bf Etape 3.}
On note $\Sigma _1$ l'\'eventail obtenu par subdivision
\'etoil\'ee de $\Sigma ''$ par rapport \`a $\rho$.
Comme $\tau$ est cod\'efinie par rapport
\`a tous les c\^ones 
maximaux de $\Sigma ''$ qui la contiennent, 
on a $\pimp (\Sigma _1) \leq \pimp (\Sigma '')$
d'apr\`es la proposition \ref{codefinie},
et puisque 
$\tau$ est contenue dans une face $\gamma$ de $\pi$-multiplicit\'e maximale, 
on a en fait $\pimp (\Sigma _1) < \pimp (\Sigma '')$,
ce qui termine la d\'emonstration 
du th\'eor\`eme de $\pi$-d\'esingularisation.
\finpreuve

\medskip

{\bf D\'emonstration de la proposition \ref{prop finale}.}
Si $\dim (\sigma) \leq 2$, il n'y a rien \`a faire~: 
$\tau$ est cod\'efinie par rapport \`a tous les c\^ones 
maximaux de $\overline{\st (\sigma)}$ qui la contiennent.
On suppose dor\'enavant que $\dim (\sigma)> 2$.

Si $\pimul (\sigma) =1$, remarquons que si $\eta$
est un c\^one maximal contenant $\sigma$, alors
toutes les faces de codimension un
$\pi$-ind\'ependantes de $\eta$
ont m\^eme $\pi$-multiplicit\'e.
Soit $\overline{\st (\sigma)}^+$
la subdivision \'etoil\'ee positive 
de $\overline{\st (\sigma)}$ par rapport \`a $\sigma$.
Alors $\pimp (\overline{\st (\sigma)}^+) \leq 
\pimp (\overline{\st (\sigma)})$ par le lemme \ref{lemme fondamental}
et la face $\tau$ est 
cod\'efinie par rapport \`a tous les c\^ones 
maximaux de $\overline{\st (\sigma)}^+$ qui la contiennent
d'apr\`es le lemme \ref{lemme codefini}.

\medskip

La d\'emonstration de la proposition \ref{prop finale}
se fait alors par r\'ecurrence sur 
$\pimul (\sigma)$~: subdivisons une premi\`ere fois
$\overline{\st (\sigma)}$ \`a l'aide de la proposition \ref{proposition
fondamentale}.
On a $\pimp (\overline{\st (\sigma)}') \leq 
\pimp (\overline{\st (\sigma)})$ par le lemme \ref{lemme fondamental}.

\smallskip

Si le cas A) de la proposition \ref{proposition
fondamentale} se produit, les circuits contenus
dans $\overline{\st (\sigma)}'$ sont tous de 
$\pi$-multiplicit\'e strictement inf\'erieure \`a celle de $\sigma$
et on conclut par l'hypoth\`ese de r\'ecurrence.

\smallskip

Si le cas B) de la proposition \ref{proposition
fondamentale} se produit,
seul le circuit $\bar{\sigma}$ contenu  
dans $\kappa '$ pose un probl\`eme {\em a priori}
car on ne peut pas lui appliquer l'hypoth\`ese de r\'ecurrence. 
Mais si $\tau$ est face d'un c\^one maximal contenant
$\kappa '$, alors $\tau \cap \sigma$ est incluse dans $\gamma '$.
En effet~: \'ecrivons 
$\sigma = \langle\rho_1,\ldots,\rho_k\rangle$ et quitte
\`a renum\'eroter les $\rho_i$, on peut supposer que
$\tau = \langle\rho_l,\ldots,\rho_k,\ldots,\rho_m\rangle$ pour
certains indices $l$ et $m$ tels que $2\leq l\leq k$ et 
$m\geq k$ (si $\tau \cap \sigma = \emptyset$, alors $\tau$
est cod\'efinie par rapport \`a tous les c\^ones 
maximaux de $\overline{\st (\sigma)}$ qui la contiennent).    
Le c\^one $\kappa '$ est de la forme 
$\langle\rho,\rho _1, \ldots,\check{\rho _i} ,\ldots, \rho _k\rangle$
pour un certain $i$ tel que $1\leq i\leq k$
et alors $\gamma ' = 
\langle\rho _1, \ldots,\check{\rho _i} ,\ldots, \rho _k\rangle$.
Comme $\tau$ est face d'un c\^one maximal contenant
$\kappa '$, on a $i \leq l-1$, c'est-\`a-dire 
$\gamma ' = \langle\rho _1, \ldots,\check{\rho _i},\ldots,\rho_l,\ldots, \rho _k\rangle$
donc $$\tau \cap \sigma = \langle\rho_l,\ldots,\rho_k\rangle \subset 
\gamma ' = \langle\rho _1, \ldots,\check{\rho _i},\ldots,\rho_l,\ldots, \rho _k\rangle.$$ 
Finalement, comme $\gamma '$ est cod\'efinie par rapport 
\`a $\kappa '$, $\tau$ est cod\'efinie par rapport 
\`a tout c\^one maximal contenant $\bar{\sigma}$.\finpreuve

\end{document}

%% file: fig1.pstex_t
\begin{picture}(0,0)%
\epsfig{file=fig1.pstex}%
\end{picture}%
\setlength{\unitlength}{0.00050000in}%
\begingroup\makeatletter\ifx\SetFigFont\undefined
\def\x#1#2#3#4#5#6#7\relax{\def\x{#1#2#3#4#5#6}}%
\expandafter\x\fmtname xxxxxx\relax \def\y{splain}%
\ifx\x\y   
\gdef\SetFigFont#1#2#3{%
  \ifnum #1<17\tiny\else \ifnum #1<20\small\else
  \ifnum #1<24\normalsize\else \ifnum #1<29\large\else
  \ifnum #1<34\Large\else \ifnum #1<41\LARGE\else
     \huge\fi\fi\fi\fi\fi\fi
  \csname #3\endcsname}%
\else
\gdef\SetFigFont#1#2#3{\begingroup
  \count@#1\relax \ifnum 25<\count@\count@25\fi
  \def\x{\endgroup\@setsize\SetFigFont{#2pt}}%
  \expandafter\x
    \csname \romannumeral\the\count@ pt\expandafter\endcsname
    \csname @\romannumeral\the\count@ pt\endcsname
  \csname #3\endcsname}%
\fi
\fi\endgroup
\begin{picture}(7747,3944)(879,-5483)
\put(3526,-1861){\makebox(0,0)[lb]{\smash{\SetFigFont{10}{12.0}{rm}$X_2 \times \{ \infty\}$}}}
\put(1801,-2011){\makebox(0,0)[lb]{\smash{\SetFigFont{10}{12.0}{rm}$p$}}}
\put(1801,-5011){\makebox(0,0)[lb]{\smash{\SetFigFont{10}{12.0}{rm}$p$}}}
\put(3526,-4861){\makebox(0,0)[lb]{\smash{\SetFigFont{10}{12.0}{rm}$X_2\times \{0\}$}}}
\put(7126,-4486){\makebox(0,0)[lb]{\smash{\SetFigFont{10}{12.0}{rm}$\bullet$}}}
\put(7126,-2086){\makebox(0,0)[lb]{\smash{\SetFigFont{10}{12.0}{rm}$\bullet$}}}
\put(2026,-5086){\makebox(0,0)[lb]{\smash{\SetFigFont{10}{12.0}{rm}$\bullet$}}}
\put(2026,-2011){\makebox(0,0)[lb]{\smash{\SetFigFont{10}{12.0}{rm}$\bullet$}}}
\put(8626,-1786){\makebox(0,0)[lb]{\smash{\SetFigFont{10}{12.0}{rm}$X_2\times\{\infty\}$}}}
\put(7426,-4486){\makebox(0,0)[lb]{\smash{\SetFigFont{10}{12.0}{rm}$E$}}}
\put(8401,-5161){\makebox(0,0)[lb]{\smash{\SetFigFont{10}{12.0}{rm}$X_1$}}}
\put(6901,-2086){\makebox(0,0)[lb]{\smash{\SetFigFont{10}{12.0}{rm}$p$}}}
\end{picture}

%% file: fig3.pstex_t
\begin{picture}(0,0)%
\epsfig{file=fig3.pstex}%
\end{picture}%
\setlength{\unitlength}{0.00050000in}%
\begingroup\makeatletter\ifx\SetFigFont\undefined
\def\x#1#2#3#4#5#6#7\relax{\def\x{#1#2#3#4#5#6}}%
\expandafter\x\fmtname xxxxxx\relax \def\y{splain}%
\ifx\x\y   
\gdef\SetFigFont#1#2#3{%
  \ifnum #1<17\tiny\else \ifnum #1<20\small\else
  \ifnum #1<24\normalsize\else \ifnum #1<29\large\else
  \ifnum #1<34\Large\else \ifnum #1<41\LARGE\else
     \huge\fi\fi\fi\fi\fi\fi
  \csname #3\endcsname}%
\else
\gdef\SetFigFont#1#2#3{\begingroup
  \count@#1\relax \ifnum 25<\count@\count@25\fi
  \def\x{\endgroup\@setsize\SetFigFont{#2pt}}%
  \expandafter\x
    \csname \romannumeral\the\count@ pt\expandafter\endcsname
    \csname @\romannumeral\the\count@ pt\endcsname
  \csname #3\endcsname}%
\fi
\fi\endgroup
\begin{picture}(11175,6912)(301,-8131)
\put(4951,-3886){\makebox(0,0)[lb]{\smash{\SetFigFont{10}{12.0}{rm}$\pi(\rho_2)$}}}
\put(5026,-5161){\makebox(0,0)[lb]{\smash{\SetFigFont{10}{12.0}{rm}$\pi(\rho_1)$}}}
\put(6901,-4186){\makebox(0,0)[lb]{\smash{\SetFigFont{10}{12.0}{rm}$\pi(\rho_4)$}}}
\put(7126,-1411){\makebox(0,0)[lb]{\smash{\SetFigFont{10}{12.0}{rm}$\pi(\rho_2)$}}}
\put(9001,-1711){\makebox(0,0)[lb]{\smash{\SetFigFont{10}{12.0}{rm}$\pi(\rho_4)$}}}
\put(9376,-3061){\makebox(0,0)[lb]{\smash{\SetFigFont{10}{12.0}{rm}$\pi(\rho_3)$}}}
\put(8326,-2311){\makebox(0,0)[lb]{\smash{\SetFigFont{10}{12.0}{rm}$v$}}}
\put(9826,-3586){\makebox(0,0)[lb]{\smash{\SetFigFont{10}{12.0}{rm}$\psi_+$}}}
\put(9676,-6361){\makebox(0,0)[lb]{\smash{\SetFigFont{10}{12.0}{rm}$\varphi_+$}}}
\put(1576,-5986){\makebox(0,0)[lb]{\smash{\SetFigFont{10}{12.0}{rm}$\rho_1$}}}
\put(1051,-4711){\makebox(0,0)[lb]{\smash{\SetFigFont{10}{12.0}{rm}$\rho_2$}}}
\put(3076,-5236){\makebox(0,0)[lb]{\smash{\SetFigFont{10}{12.0}{rm}$\rho_3$}}}
\put(2401,-3811){\makebox(0,0)[lb]{\smash{\SetFigFont{10}{12.0}{rm}$\rho_4$}}}
\put(7201,-2836){\makebox(0,0)[lb]{\smash{\SetFigFont{10}{12.0}{rm}$\pi(\rho_1)$}}}
\put(6376,-3211){\makebox(0,0)[lb]{\smash{\SetFigFont{10}{12.0}{rm}$\psi_-$}}}
\put(301,-4786){\makebox(0,0)[lb]{\smash{\SetFigFont{10}{12.0}{rm}$\pi$}}}
\put(7951,-4411){\makebox(0,0)[lb]{\smash{\SetFigFont{10}{12.0}{rm}$\varphi$}}}
\put(5926,-6511){\makebox(0,0)[lb]{\smash{\SetFigFont{10}{12.0}{rm}$\varphi_-$}}}
\put(7951,-3211){\makebox(0,0)[lb]{\smash{\SetFigFont{10}{12.0}{rm}$X_v$}}}
\put(1126,-6511){\makebox(0,0)[lb]{\smash{\SetFigFont{10}{12.0}{rm}$\sigma = \tau \subset N^+$}}}
\put(7651,-8086){\makebox(0,0)[lb]{\smash{\SetFigFont{10}{12.0}{rm}$X_{\tau}/\!/\KK^*$}}}
\put(7276,-6286){\makebox(0,0)[lb]{\smash{\SetFigFont{10}{12.0}{rm}$\pi(\rho_2)$}}}
\put(9151,-6811){\makebox(0,0)[lb]{\smash{\SetFigFont{10}{12.0}{rm}$\pi(\rho_4)$}}}
\put(9376,-7786){\makebox(0,0)[lb]{\smash{\SetFigFont{10}{12.0}{rm}$\pi(\rho_3)$}}}
\put(11476,-5386){\makebox(0,0)[lb]{\smash{\SetFigFont{10}{12.0}{rm}$\pi(\rho_3)$}}}
\put(7276,-7561){\makebox(0,0)[lb]{\smash{\SetFigFont{10}{12.0}{rm}$\pi(\rho_1)$}}}
\put(8926,-3811){\makebox(0,0)[lb]{\smash{\SetFigFont{10}{12.0}{rm}$\pi(\rho_2)$}}}
\put(11251,-4186){\makebox(0,0)[lb]{\smash{\SetFigFont{10}{12.0}{rm}$\pi(\rho_4)$}}}
\put(9301,-5236){\makebox(0,0)[lb]{\smash{\SetFigFont{10}{12.0}{rm}$\pi(\rho_1)$}}}
\put(8851,-5836){\makebox(0,0)[lb]{\smash{\SetFigFont{10}{12.0}{rm}$X(\pi(\partial_+(\sigma)))=(X_{\tau})_+/\KK^*$}}}
\put(7276,-5386){\makebox(0,0)[lb]{\smash{\SetFigFont{10}{12.0}{rm}$\pi(\rho_3)$}}}
\put(4426,-5761){\makebox(0,0)[lb]{\smash{\SetFigFont{10}{12.0}{rm}$X(\pi(\partial_-(\sigma)))=(X_{\tau})_-/\KK^*$}}}
\end{picture}

%% file: fig2.pstex_t
\begin{picture}(0,0)%
\epsfig{file=fig2.pstex}%
\end{picture}%
\setlength{\unitlength}{0.00058300in}%
\begingroup\makeatletter\ifx\SetFigFont\undefined
\def\x#1#2#3#4#5#6#7\relax{\def\x{#1#2#3#4#5#6}}%
\expandafter\x\fmtname xxxxxx\relax \def\y{splain}%
\ifx\x\y   
\gdef\SetFigFont#1#2#3{%
  \ifnum #1<17\tiny\else \ifnum #1<20\small\else
  \ifnum #1<24\normalsize\else \ifnum #1<29\large\else
  \ifnum #1<34\Large\else \ifnum #1<41\LARGE\else
     \huge\fi\fi\fi\fi\fi\fi
  \csname #3\endcsname}%
\else
\gdef\SetFigFont#1#2#3{\begingroup
  \count@#1\relax \ifnum 25<\count@\count@25\fi
  \def\x{\endgroup\@setsize\SetFigFont{#2pt}}%
  \expandafter\x
    \csname \romannumeral\the\count@ pt\expandafter\endcsname
    \csname @\romannumeral\the\count@ pt\endcsname
  \csname #3\endcsname}%
\fi
\fi\endgroup
\begin{picture}(9624,9126)(889,-8875)
\put(3001,-961){\makebox(0,0)[lb]{\smash{\SetFigFont{11}{13.2}{rm}$\eta$}}}
\put(2401,-661){\makebox(0,0)[lb]{\smash{\SetFigFont{11}{13.2}{rm}$\bullet$}}}
\put(5926,-1336){\makebox(0,0)[lb]{\smash{\SetFigFont{11}{13.2}{rm}Etape 1}}}
\put(8701,-1261){\makebox(0,0)[lb]{\smash{\SetFigFont{11}{13.2}{rm}$\bullet$}}}
\put(8701,-1486){\makebox(0,0)[lb]{\smash{\SetFigFont{11}{13.2}{rm}$\rho$}}}
\put(2401,-811){\makebox(0,0)[lb]{\smash{\SetFigFont{11}{13.2}{rm}$\rho_{\pm}$}}}
\put(8026,-811){\makebox(0,0)[lb]{\smash{\SetFigFont{11}{13.2}{rm}$\rho_{\pm}$}}}
\put(9151,-1036){\makebox(0,0)[lb]{\smash{\SetFigFont{11}{13.2}{rm}$\tau$}}}
\put(8551,-886){\makebox(0,0)[lb]{\smash{\SetFigFont{11}{13.2}{rm}$\nu$}}}
\put(1051,-1861){\makebox(0,0)[lb]{\smash{\SetFigFont{11}{13.2}{rm}$\pi$}}}
\put(5851,-6511){\makebox(0,0)[lb]{\smash{\SetFigFont{11}{13.2}{rm}Etape 3}}}
\put(8701,-6361){\makebox(0,0)[lb]{\smash{\SetFigFont{11}{13.2}{rm}$\bullet$}}}
\put(9226,-6061){\makebox(0,0)[lb]{\smash{\SetFigFont{11}{13.2}{rm}$\tau$}}}
\put(1726,-2536){\makebox(0,0)[lb]{\smash{\SetFigFont{11}{13.2}{rm}$\Sigma$}}}
\put(7051,-2386){\makebox(0,0)[lb]{\smash{\SetFigFont{11}{13.2}{rm}$\Sigma'$}}}
\put(2401,-3211){\makebox(0,0)[lb]{\smash{\SetFigFont{11}{13.2}{rm}$\pimp (\eta)$ est maximal }}}
\put(9151,-4261){\makebox(0,0)[lb]{\smash{\SetFigFont{11}{13.2}{rm}Etape 2}}}
\put(9376,-1711){\makebox(0,0)[lb]{\smash{\SetFigFont{11}{13.2}{rm}$\sigma'$}}}
\put(7051,-7336){\makebox(0,0)[lb]{\smash{\SetFigFont{11}{13.2}{rm}$\Sigma''$}}}
\put(8476,-6136){\makebox(0,0)[lb]{\smash{\SetFigFont{11}{13.2}{rm}$\nu$}}}
\put(3226,-6586){\makebox(0,0)[lb]{\smash{\SetFigFont{11}{13.2}{rm}$\rho$}}}
\put(8701,-6511){\makebox(0,0)[lb]{\smash{\SetFigFont{11}{13.2}{rm}$\rho$}}}
\put(7351,-3136){\makebox(0,0)[lb]{\smash{\SetFigFont{11}{13.2}{rm}$\pimp (\nu)$ est maximal et $\tau$ est  }}}
\put(7351,-3436){\makebox(0,0)[lb]{\smash{\SetFigFont{11}{13.2}{rm}cod\'efinie par rapport \`a $\nu$}}}
\put(6751,-8236){\makebox(0,0)[lb]{\smash{\SetFigFont{11}{13.2}{rm}$\pimp(\nu)$ est maximal et $\tau$ est }}}
\put(6751,-8536){\makebox(0,0)[lb]{\smash{\SetFigFont{11}{13.2}{rm}cod\'efinie par rapport \`a tout c\^one }}}
\put(6751,-8836){\makebox(0,0)[lb]{\smash{\SetFigFont{11}{13.2}{rm}maximal contenant $\tau$}}}
\put(1201,-7186){\makebox(0,0)[lb]{\smash{\SetFigFont{11}{13.2}{rm}$\Sigma_1$}}}
\put(1801,-8311){\makebox(0,0)[lb]{\smash{\SetFigFont{11}{13.2}{rm}$\pimp (\Sigma_1) < \pimp(\Sigma)$}}}
\end{picture}